\documentstyle[12pt]{article}

\begin{document}
\def\comment#1{\marginpar{{\scriptsize #1}}}
\def\framew#1{\fbox{#1}}
\def\framep#1{\noindent \fbox{\vbox{#1}}}
\def\framef#1{\fbox{\vbox{ #1 }}}
\def\be{\begin{equation}}
\def\bea{\begin{eqnarray}}
\def\ee{\end{equation}}
\def\eea{\end{eqnarray}}
\newcommand{\sect}[1]{\setcounter{equation}{0}\section{#1}}
\newcommand{\subsect}[1]{\subsection{#1}}
\newcommand{\subsubsect}[1]{\subsubsection{#1}}
\renewcommand{\theequation}{\thesection.\arabic{equation}}
\renewcommand{\thefootnote}{\fnsymbol{footnote}}
\def\fn{\footnote}
\footskip 1.0cm
\def\sxn#1{\bigskip\medskip \sect{#1} \smallskip}
\def\subsxn#1{\medskip \subsect{#1} \smallskip}
\def\subsubsxn#1{\medskip \subsubsect{#1} \smallskip}
\newtheorem{proposition}{Proposition}[section]  
\def\bprop{\medskip\begin{proposition}~~~\smallskip \it}
\def\eprop{\end{proposition}\bigskip}
\def\proof{\bigskip \noindent {\it Proof.} \ }
\newtheorem{naming}{Definition}[section]   
\def\bnam{\medskip\begin{naming}~~~\smallskip \it}
\def\enam{end{naming}\bigskip}
\newtheorem{example}{Example}[section]   
\def\bexam{\medskip\begin{example} ~~\\ \rm}
\def\eexam{\end{example}\bigskip}
\bibliographystyle{unsrt}
\def\br{\begin{thebibliography}{abc}}
\def\er{\end{thebibliography}}
\def\rf{\bibitem}

%
%
\def\pinco{| \! / \! |}
\def\cstars{$C^*$-algebras }
\def\cstar{$C^*$-algebra }
\def\unit{I\!\!I}
\def\norm#1{\parallel #1 \parallel}
\def\abs#1{\left| #1\right|}
\def\ha{\widehat{\cal A}}
\def\hc{\widehat{\cal C}}
\def\prim{$Prim \ca$ }
\def\bar#1{\overline{#1}}
\def\what{\widehat}
\def\wtilde{\wtilde}
\def\bra#1{\left\langle #1\right|}
\def\ket#1{\left| #1\right\rangle}
\def\EV#1#2{\left\langle #1\vert #2\right\rangle}
\def\VEV#1{\left\langle #1\right\rangle}
%
%
\def\sdp{\hbox{ \raisebox{.25ex}{\tiny $|$}\hspace{.2ex}{$\!\!\times $}} }
\def\pa{\partial}
\def\del{\nabla}
\def\a{\alpha}
\def\b{\beta}
\def\c{\raisebox{.4ex}{$\chi$}}
\def\d{\delta}
\def\e{\epsilon}
\def\f{\phi}
\def\g{\gamma}
\def\h{\eta}
\def\i{\iota}
\def\j{\psi}
\def\k{\kappa}
\def\l{\lambda}
\def\m{\mu}
\def\n{\nu}
\def\o{\omega}
\def\p{\pi}
\def\q{\theta}
\def\r{\rho}
\def\s{\sigma}
\def\t{\tau}
\def\u{\upsilon}
\def\x{\xi}
\def\z{\zeta}
\def\D{\Delta}
\def\F{\Phi}
\def\G{\Gamma}
\def\J{\Psi}
\def\L{\Lambda}
\def\O{\Omega}
\def\P{\Pi}
\def\Q{\Theta}
\def\S{\Sigma}
\def\U{\Upsilon}
\def\X{\Xi}
\def\Z{\Zeta}
\def\ca{{\cal A}}
\def\cb{{\cal B}}
\def\cc{{\cal C}}
\def\cd{{\cal D}}
\def\ce{{\cal E}}
\def\cf{{\cal F}}
\def\cg{{\cal G}}
\def\ch{{\cal H}}
\def\ci{{\cal I}}
\def\cj{{\cal J}}
\def\ck{{\cal K}}
\def\cl{{\cal L}}
\def\cm{{\cal M}}
\def\cn{{\cal N}}
\def\co{{\cal O}}
\def\cp{{\cal P}}
\def\ct{{\cal T}}
\def\cu{{\cal U}}
\def\cv{{\cal V}}
\def\cw{{\cal W}}
\def\cx{{\cal X}}
\def\cy{{\cal }}
\def\cz{{\cal Z}}
%
%
%
\def\inbar{\,\vrule height1.5ex width.4pt depth0pt}
\def\IC{\relax\,\hbox{$\inbar\kern-.3em{\rm C}$}}
\def\ID{\relax{\rm I\kern-.18em D}}
\def\IF{\relax{\rm I\kern-.18em F}}
\def\IH{\relax{\rm I\kern-.18em H}}
\def\II{\relax{\rm I\kern-.17em I}}
\def\I1{\relax{\rm 1\kern-.28em l}}
\def\IM{\relax{\rm I\kern-.18em M}}
\def\IN{\relax{\rm I\kern-.18em N}}
\def\IP{\relax{\rm I\kern-.18em P}}
\def\IQ{\relax\,\hbox{$\inbar\kern-.3em{\rm Q}$}}
\def\IZ{\relax\,\hbox{$\inbar\kern-.3em{\rm Z}$}}
\def\IR{\relax{\rm I\kern-.18em R}}
\font\cmss=cmss10 \font\cmsss=cmss10 at 7pt
\def\Z{\relax\ifmmode\mathchoice
{\hbox{\cmss Z\kern-.4em Z}}{\hbox{\cmss Z\kern-.4em Z}}
{\lower.9pt\hbox{\cmsss Z\kern-.4em Z}}
{\lower1.2pt\hbox{\cmsss Z\kern-.4em Z}}\else{\cmss Z\kern-.4emZ}\fi}
%
%
\def\bc{{\bf C}}
\def\br{{\bf R}}
\def\bz{{\bf Z}}
\def\bn{{\bf N}}
\def\bm{{\bf M}}
\def\Up{\Uparrow}
\def\up{\uparrow}
\def\Dn{\Downarrow}
\def\dn{\downarrow}
\def\Ra{\Rightarrow}
\def\ra{\rightarrow}
\def\La{\Leftarrow}
\def\la{\leftarrow}
\def\iff{\Leftrightarrow}
\thispagestyle{empty}
\setcounter{page}{0}
\hfill \today

\vspace{1cm}

\begin{center}{\Large $K$-THEORY OF NONCOMMUTATIVE LATTICES}
\end{center}
\vspace{1cm}
\centerline {Elisa Ercolessi$^1$,
             Giovanni Landi$^{2,3}$,
	     Paulo Teotonio-Sobrinho$^{2,4}$}
\vspace{1cm}
\centerline {\it $^1$  Dipartimento di Fisica,
Universit\`a di Bologna, and INFM}
\centerline{\it Via Irnerio 46, I-40126, Bologna, Italy.}
\vspace{2mm}
\centerline {\it $^2$ The E. Schr\"odinger International Institute for
Mathematical Physics,}
\centerline{\it Pasteurgasse 6/7, A-1090 Wien, Austria.}
\vspace{2mm}
\centerline{\it $^3$ Dipartimento di Scienze Matematiche,
Universit\`a di Trieste,}
\centerline{\it P.le Europa 1, I-34127, Trieste, Italy.}
\centerline {\it and ~ INFN, Sezione di Napoli, Napoli, Italy.}
\vspace{2mm}
\centerline {\it $^4$ Dept.\ of Physics, Univ.\ of  Illinois at Chicago,}
\centerline{\it 60607-7059, Chicago, IL, USA.}
\centerline {\it and Universidade de Sao Paulo, Instituto de Fisica - DFMA,}
\centerline{\it Caixa Postal 66318, 05389-970, Sao Paulo, SP,
Brasil \footnote{Permanent address}.} 
\vspace{.5cm}
\begin{abstract}

Noncommutative lattices have been recently used as finite topological 
approximations in quantum physical models. As a first step in the 
construction of bundles and characteristic classes over such 
noncommutative spaces, we shall study their $K$-theory. We shall do it
algebraically, by studying the algebraic $K$-theory of the associated 
algebras of `continuous functions' which turn out to be noncommutative
approximately finite dimensional (AF) $C^*$-algebra . We also work out 
several examples.
\end{abstract}

\vfill\eject
\setcounter{page}{1}

\renewcommand{\thefootnote}{\arabic{footnote}}
\setcounter{footnote}{0}

\sxn{Introduction}\label{se:int}

Topological lattices, namely finite topological spaces, were introduced in
\cite{So} (with the name posets or partially ordered sets) as finite
topological 
approximations to `continuum' spaces.  Their ability to capture some of the
topological information of the space they approximate has been the motivation
for their use in quantum theories \cite{BBET,lattice,BBLLT}. The idea is to
construct alternative lattice theories which are able to describe also 
some topologically non trivial configurations. An example of a promising 
result in
this direction is the construction of $\q$-states for particles on the poset
approximations to a circle. In a suitable limit these states give the usual
$\q$-states for particles on the circle \cite{BBET,lattice}. Also, there is
some work on gauge field theories \cite{BBLLT}.

In \cite{lattice} it was observed that a poset $P$ is truly a noncommutative
space (in fact a noncommutative lattice) since it can be described as the 
structure space (space of irreducible representations) of a noncommutative
\cstar $\ca $. Therefore, such an algebra is the analogue of the
commutative algebra of complex valued continuous functions defined on any
Hausdorff topological space. The algebra $\ca $ can be thought of as an algebra
of operator valued functions over $P$.
The use of these algebras leads to noncommutative geometry \cite{Co} as
the natural tool to construct ``geometric" structures on noncommutative
lattices. It turns out that the algebras which are relevant are approximately
finite dimensional (AF) postliminal algebras. 
There are in general several algebras which are associated with the same 
space $P$ and several ways of constructing any such an algebra
\cite{ELTfunctions}.  In this paper we will use a
diagrammatic method due to Bratteli \cite{Br1}. Although this method gives only
one $\ca$ among all possibilities, it will be enough for our purposes.

In this paper, as a first step in the construction of bundles and characteristi
classes over noncommutative lattices, we shall study the $K$-theory of such
noncommutative lattices. We shall do it from the algebraic point of view, by
studying the algebraic
$K$-theory of the associated algebras. Algebraically, the analogue of vector
bundles over a noncommutative lattice $P$ are projective modules of finite type
over the corresponding algebra $\ca$. The $K$-theory group of $\ca$  classifies
(stable) equivalence classes of such projective modules. For the class of
algebras we are interested in, Bratteli diagrams allow to explicitly construct
these groups.

The paper is organized as follows. In Section~\ref{se:top} we shall recall the
construction of topological lattices as approximation spaces of
`continuum' topological spaces. Section~\ref{se:ncl} is devoted to the
description of the noncommutative algebras associated with the noncommutative
lattices, and the description of the Bratteli diagrams. Several examples are
worked out. Section~\ref{se:pbk} treats the notion of projective modules and
the $K$-theory. Again, several examples are presented.

\sxn{Topological Lattices}\label{se:top}

We shall briefly recall how to construct a topological lattice from
any covering 
of a `continuum' topological space $M$ while referring to
\cite{So,lattice,ELTfunctions} for details. The idea is simply to identify any
two points of $M$ which cannot be separated or distinguished by the sets in the
covering. The resulting space consists of a finite (or in general a countable) 
number of points with a nontrivial topology which maintain some of the
topological information of the starting space $M$ \cite{So}.

If we are given a covering $ \cu =\{O_\l \}$ of $M$ which is a topology for the
latter, any two points $x$ and $y$ in $M$, will be declared to be equivalent,
$x\sim y$, if every set $O_\l$ containing either point $x$ or $y$ contains the
other too,
\be
   x\sim y ~~~~{\rm if~ and~ only~ if}~~~~ x\in O_\l \iff y\in O_\l
~~~\forall~~ O_\l~ . \label{2.2}
\ee
Then, we replace $M$ by $P(M) =: M /\sim$ and endow this space with
the quotient 
topology. 
When $M$ is compact, the number of sets $O_\l$ in
$\cu$ can be taken to be finite so that $P(M)$ is an approximation to $M$ by a
finite number of points. When $M$ is not compact we take it locally compact so
that each point has a neighborhood intersected by only finitely many $O_\l$ and
$P(M)$ is a ``finitary" (countable) approximation to $M$ \cite{So}. 

In general, the space $P(M)$ is neither  Hausdorff (one cannot separate
completely any two points) nor $T_1$ (not all points are
closed). However, it is 
always a a $T_0$ space \cite{So}. This means that given any two points, there
exists at least an open set which contains only one of the points and not the
other.  

The space $P(M)$ is made a {\it partially ordered set} (or a {\it poset})
with a partial order $\preceq$ given by
\[ x\preceq y~~~ \mbox{ if every open set containing } y
\mbox{ contains also } x~.\]
The smallest open set $O_x$ containing a point $x \in P(M)$ consists
of all $y$'s 
preceding $x$: $O_x = \{y \in P(M) \; : \; y \preceq x \}$. The open sets
$\{O_x~,  x \in P(M) \}$ are a basis for the topology of $P(M)$. For spaces
with at most a countable number of points, a topology is equivalent
to a partial order. 
Any poset can be represented pictorially by a {\it Hasse diagram}
constructed by 
arranging its points at different levels and connecting them using the
following rules:
\begin{enumerate}
\item[1)] if $x\prec y$, then $x$ is at a lower level than $y$;
\item[2)] if $x\prec y$ and there is no $z$ such that $x\prec z\prec y$,
then $x$ is
at the level immediately below $y$ and these two points
are connected by a line called a link.
\end{enumerate}

Fig.~\ref{fi:cirhas} shows the Hasse diagram for the poset $P_4(S^1)$, a four
point approximation to $S^1$. A basis of open sets for the topology is given by
\be
\{x_1\}~,  ~ ~ \{x_2\}~,  ~ ~ \{x_1,x_2,x_3\}~,
~ ~ \{x_1,x_2, x_4\}~. \label{2.5}
\ee
\begin{figure}[htb]
\begin{center}
\begin{picture}(230,60)(-110,-30)
\put(-30,30){\circle*{4}}
\put(30,30){\circle*{4}}
\put(-30,-30){\circle*{4}}
\put(30,-30){\circle*{4}}
\put(-30,30){\line(0,-1){60}}
\put(30,30){\line(0,-1){60}}
\put(-30,30){\line(1,-1){60}}
\put(30,30){\line(-1,-1){60}}
\put(-45,29){$x_3$}
\put(35,29){$x_4$}
\put(-45,-33){$x_1$}
\put(35,-33){$x_2$}
\end{picture}
\end{center}
\caption{\label{fi:cirhas} 
\protect{\footnotesize The Hasse diagram for $P_4(S^1)$.  }}
\end{figure}
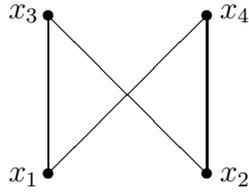\bigskip

Fig.~\ref{fi:sphpos} shows  the Hasse diagram for the poset $P_6(S^2)$, a six
points approximation to $S^2$. A basis of open sets for the topology
is given by  
\bea
&& \{x_1\}~,  ~ ~ \{x_2\}~,  ~ ~ \{x_1,x_2,x_3\}~,
~ ~ \{x_1,x_2,x_4\}~, \nonumber \\
&& \{x_1,x_2,x_3,x_4,x_5\}~, ~ ~ \{x_1,x_2,x_3,x_4,x_6\}~.
~ ~ \label{2.6}
\eea

As alluded to before, posets maintain some of the topological information of
the spaces they approximate. This is shown, for instance, by their ability
to reproduce homotopy groups. For example, $\p_1(P_N(S^1)) = \IZ$ whenever 
$N \geq 4$ \cite{So}. Similarly, $\p_1(P_6(S^2)) =\{0\}$ and 
$\p_2(P_6(S^2)) = \IZ$. It is this fact that makes posets relevant for
quantum physical models such as $\q$-states \cite{BBET,lattice}.
Finally, we mention that the topological space being approximated
can be recovered by considering a sequence of finer and finer coverings, the
appropriated framework being that of projective systems of topological spaces.
We refer to \cite{So,pangs} for details. 
 
\begin{figure}[htb]
\begin{center}
\begin{picture}(200,120)(-110,-30)
\put(-30,90){\circle*{4}}
\put(30,90){\circle*{4}}
\put(-30,30){\circle*{4}}
\put(30,30){\circle*{4}}
\put(-30,-30){\circle*{4}}
\put(30,-30){\circle*{4}}
\put(-30,90){\line(0,-1){60}}
\put(30,90){\line(0,-1){60}}
\put(-30,90){\line(1,-1){60}}
\put(30,90){\line(-1,-1){60}}
\put(-30,30){\line(0,-1){60}}
\put(30,30){\line(0,-1){60}}
\put(-30,30){\line(1,-1){60}}
\put(30,30){\line(-1,-1){60}}
\put(-45,90){$x_5$}
\put(35,90){$x_6$}
\put(-45,29){$x_3$}
\put(35,29){$x_4$}
\put(-45,-33){$x_1$}
\put(35,-33){$x_2$}
\end{picture}
\end{center}
\caption{\label{fi:sphpos} 
\protect{\footnotesize The Hasse diagram for the poset $P_6(S^2)$.   }}
\end{figure}
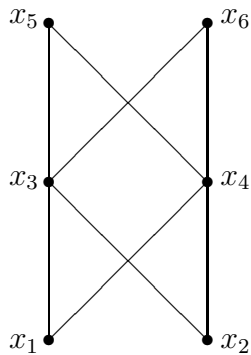\bigskip

\sxn{Noncommutative Lattices}\label{se:ncl}

Associated with any poset $P$ there is a noncommutative \cstar $\ca$ of
operator valued continuous functions on $P$. The latter can be
identified with the space $\ha = Prim\ca$ of primitive ideals of $\ca$, 
an ideal of $\ca$ being called primitive if it is the kernel of an irreducible
representation of $\ca$. On $\ha$ there is a partial order defined 
by inclusion of ideals: given any two ideals $I_1, I_2 \in \ha$, $I_1 \preceq
I_2$, if and only if $I_1 \subseteq I_2$. The space $P$ being at most
countable, the resulting topology coincides with the Jacobson topology.
Thus any poset is a noncommutative space (in fact a {\it noncommutative
lattice}) \cite{Co}. 
 
It turns out that the algebras one is working with are particularly simple.
They are approximately finite dimensional (AF) postliminal algebras. 
An AF algebra can be approximated in norm by direct sums of finite dimensional
matrix algebras.  Postliminal algebras have the remarkable property that
their irreducible (unitary) representations are completely characterized
by their kernels \cite{Di,Mu}.

In the rest of this section we shall briefly describe a diagrammatic
representation of AF algebras due to Bratteli \cite{Br1} which is very useful
for the study of the $K$-theory of the
associated posets. We refer to \cite{ELTfunctions} for a detailed treatment.

\subsxn{The Algebra of a Poset and Bratteli Diagrams}\label{se:bra}

A \cstar $\ca$ is said to be {\it{approximately finite dimensional}}
(AF) if 
there exists an increasing sequence
\be
\ca_0 ~{\buildrel I_0 \over \hookrightarrow}~ \ca_1
      ~{\buildrel I_1 \over \hookrightarrow}~ \ca_2
      ~{\buildrel I_2 \over \hookrightarrow}~ \cdots
      ~{\buildrel I_{n-1} \over \hookrightarrow}~ \ca_n
      ~{\buildrel I_n \over \hookrightarrow} \cdots
\label{af}
\ee
of finite dimensional $C^*$-subalgebras $\ca_n$ of $\ca$, with injective
$^*$-homomorphisms $I_n$, such that $\ca$ is the norm closure of $\bigcup_n
\ca_n$.  Here the maps $I_n$ are injective
$^*$-homomorphisms. Elements of $\ca$ 
are then coherent sequences of the form,
\be
a = (a_n)_{n \in \IN}~, ~a_n \in \ca_n ~~{\rm s.t.}~~ \exists  N_0 ~,
~a_{n+1} =  I_n(a_n)~, \forall ~n>N_0.
\ee
The norm of any such element is given by
\be
||(a_n)_{n \in \IN}|| = \lim_{n \ra \infty} ||a_n||_{\ca_n}~. \label{norm}
\ee
Since the maps $I_n$ are injective, the expression (\ref{norm}) gives
directly a true norm and not simply a seminorm and there is no need to
quotient out the zero norm elements \cite{W-O}.

There is a very useful diagrammatic representation of the AF algebra (\ref{af})
due to Bratteli \cite{Br1}. Each subalgebra $\ca_n$ is a direct sum of matrix
algebras
\be
\ca_n = \bigoplus_{k=1}^{n_n} M^{(n)}(d_k,\IC)~,
\ee
with $M^{(n)}(d_k,\IC)$ the algebra of $d_k \times d_k$ matrices with complex
coefficients. Given any two such algebras 
$\ca_1 = \bigoplus_{j=1}^{n_1} \IM(d_j^{(1)},\IC)$ and
$\ca_2 = \bigoplus_{k=1}^{n_2} \IM(d_k^{(2)},\IC)$ with an embedding 
$\ca_1 \hookrightarrow \ca_2$, one can always choose suitable bases in 
$\ca_1$ and $\ca_2$ in such a manner to identify $\ca_1$ with a subalgebra of
$\ca_2$ of the following form
\be
\ca_1 \simeq \bigoplus_{k=1}^{n_2} \left( \bigoplus_{j=1}^{n_1} N_{kj}
\IM(d_j^{(1)},\IC) \right)  \; .
\ee
Here, with $p,q$ any two nonnegative integers, the symbol $p\IM(q,\IC)$ 
stands for $\IM(q,\IC) \otimes{\IC}\II_p$. The nonnegative
integers
$N_{kj}$ satisfies the condition $\sum_{j=1}^{n_1} N_{kj} d^{(1)}_j =
d^{(2)}_k$.
One says that the algebra $\IM(d_j^{(1)},\IC)$ is {\it partially embedded} in
$\IM(d_k^{(2)},\IC)$ with {\it multiplicity} $N_{kj}$. 
One represents the embedding $\ca_1 \hookrightarrow \ca_2$
by means of a {\it diagram} (the Bratteli diagram), which is constructed out of
the dimensions $d_j$ ($j=1,\ldots,n_1$) and $d_k$ ($k=1,\ldots,n_2$) of the
diagonal blocks of the two algebras and the numbers $N_{kj}$ that describe the
embedding. To construct the diagram one draws two horizontal rows of
vertices, the top (bottom) one representing $\ca_1$ ($\ca_2$) and consisting of
$n_1$ ($n_2$) vertices labeled by $d_1, \ldots, d_{n_1}$
($d_1,\ldots,d_{n_2}$). 
Then for each $j=1,\ldots,n_1$ and $k=1,\ldots,n_2$, draw $N_{kj}$
edges between 
$d_j$ and $d_k$. One writes $d^{(K)}_j \searrow^{p} d^{(K+1)}_k$ to denote the
fact that $M^{(K)}(d^{(K)}_j,\IC)$ is embedded in $M^{(K+1)}(d^{(K+1)}_k,\IC)$
with multiplicity $p$. For any AF algebra $\ca$ one repeats the procedure for 
each level so obtaining a semi-infinite diagram denoted by $\cd(\ca)$ which
completely defines $\ca$ up to isomorphisms.

\bexam
Consider the subalgebra $\ca$ of the algebra $\cb(\ch)$ of bounded operators on
an infinite dimensional (separable) Hilbert space $\ch = \ch_1 \oplus \ch_2$,
given in the following manner. Let $\cp_j$ be the projection operators on
$\ch_j~, j = 1, 2$ and $\ck(\ch)$ be the algebra of compact operators on $\ch$.
The algebra is then 
\be
\ca(\vee) =
\IC\cp_1 + \ck({\ch}) + \IC\cp_2~. \label{alvee}
\ee
The use of the symbol $\ca(\bigvee)$ will be clear later. This \cstar can
be obtained as the direct limit of the following sequence of finite dimensional
algebras:
\bea
& & \ca_0 = M(1,\IC)   \nonumber \\
& & \ca_1 = M(1,\IC) \oplus M(1,\IC)  \nonumber \\
& & \ca_2 = M(1,\IC) \oplus M(2,\IC) \oplus M(1,\IC)  \nonumber \\
& & \ca_3 = M(1,\IC) \oplus M(4,\IC) \oplus M(1,\IC) \nonumber \\
& & ~~~ \vdots \nonumber \\
& & \ca_n = M(1,\IC) \oplus M(2n-2,\IC) \oplus M(1,\IC) \nonumber \\
& & ~~~ \vdots
\label{vee}
\eea
with $\ca_n$ embedded in $\ca_{n+1}$ as 
$M(1,\IC) \oplus (M(1,\IC) \oplus M(2n-2,\IC) \oplus M(1,\IC)) \oplus
M(1,\IC)$, 
\be
a_{n} = \left[
\begin{array}{ccc}
\l_1 & 0                  & 0    \\
0    & m_{(2n-2)\times(2n-2)} & 0    \\
0    & 0                  & \l_2
\end{array}
\right]~~ \ra ~
\left[
\begin{array}{ccccc}
\l_1 & 0    & 0                  & 0       & 0      \\
0    & \l_1 & 0                  & 0       & 0      \\
0    & 0    & m_{(2n-2)\times(2n-2)} & 0       & 0      \\
0    & 0    & 0                  & \l_2    & 0       \\
0    & 0    & 0                  & 0       & \l_2
\end{array}
\right]~.
\label{vee1}
\ee
\begin{figure}[htb]
\begin{center}
\begin{picture}(120,160)(0,40)
\put(30,150){\circle*{4}}
\put(30,120){\circle*{4}}
\put(30,90){\circle*{4}}
\put(30,60){\circle*{4}}
\put(60,180){\circle*{4}}
\put(60,120){\circle*{4}}
\put(60,90){\circle*{4}}
\put(60,60){\circle*{4}}
\put(90,150){\circle*{4}}
\put(90,120){\circle*{4}}
\put(90,90){\circle*{4}}
\put(90,60){\circle*{4}}
\put(30,150){\line(0,-1){30}}
\put(30,120){\line(0,-1){30}}
\put(30,90){\line(0,-1){30}}
\put(60,120){\line(0,-1){30}}
\put(60,90){\line(0,-1){30}}
\put(90,150){\line(0,-1){30}}
\put(90,120){\line(0,-1){30}}
\put(90,90){\line(0,-1){30}}
\put(30,60){\line(0,-1){10}}
\put(60,60){\line(0,-1){10}}
\put(90,60){\line(0,-1){10}}
\put(60,180){\line(1,-1){30}}
\put(30,150){\line(1,-1){30}}
\put(30,120){\line(1,-1){30}}
\put(30,90){\line(1,-1){30}}
\put(60,180){\line(-1,-1){30}}
\put(90,150){\line(-1,-1){30}}
\put(90,120){\line(-1,-1){30}}
\put(90,90){\line(-1,-1){30}}
\put(30,60){\line(1,-1){10}}
\put(90,60){\line(-1,-1){10}}
\put(54,185){{\small$1$}} 
\put(13,148){{\small$1$}} 
\put(13,118){{\small$1$}}
\put(13,88){{\small$1$}} 
\put(13,58){{\small$1$}} 
\put(93,148){{\small$1$}}
\put(93,118){{\small$1$}} 
\put(93,88){{\small$1$}} 
\put(93,58){{\small$1$}}
\put(63,113){$2$}
\put(63,83){$4$}
\put(63,53){$6$}
\put(45,30){$\vdots$}
\put(75,30){$\vdots$}
\end{picture}
\end{center}
\caption{\label{fi:veealg} 
\protect{\footnotesize The Bratteli diagram of the algebra $\ca(\bigvee)$. }}
\end{figure}\bigskip

\noindent
The corresponding Bratteli diagram is shown in Fig.~\ref{fi:veealg}.
\eexam

Out of the Bratteli diagram $\cd(\ca)$ of an AF algebra $\ca$ one can also
identify the ideals of $\ca$ and decide which ones are primitive. As we have
mentioned before, the topology is given by constructing a poset whose partial
order is provided by inclusion of ideal. In \cite{Br2} it is proved that the
(norm closed) ideals $\{\ci\}$ of $\ca$ are all and only the (norm closure of) 
sums of the form 
\be
\ci = \bigcup_{n=1}^{\infty} \oplus_{k, (n,k) \in \Lambda}
~\IM^{(n)}(d_k,\IC)~, \label{idea}  
\ee
with the subset $\Lambda \equiv \Lambda_{\ci}$ of the Bratteli diagram
satisfying 
the following two properties:
\begin{itemize}
\item[i)]
if $M^{(n)}(d_k,\IC)\in \Lambda$ and $M^{(n)}(d_k,\IC)
\searrow M^{(n+1)}(d_j,\IC)$ then  $M^{(n+1)}(d_j,\IC)$
belongs to $\Lambda$;
\item[ii)]
if all factors $M^{(n+1)}(d_j,\IC)~, j=\{1,2,\cdots,N_{n+1}\},$ in which
$M^{(n)}(d_k,\IC)$ is partially embedded belong to
$\Lambda$, then $M^{(n)}(d_k,\IC)$ belongs to $\Lambda$.
\end{itemize}
Furthermore, the ideal (\ref{idea}) is primitive if and only if the associated
subdiagram $\Lambda_{\ci}$ satisfies the additional property:
\begin{itemize}
\item[iii)]
$\forall~ n$ there exists an $M^{(m)}(d_j,\IC)$, with $m >
n$, not belonging to $\Lambda$ such that, for all $k\in\{1,2,\cdots,N_n\}$ with
$M^{(n)}(d_k,\IC)$ not in $\Lambda$, one can find a sequence
$M^{(n)}(d_k,\IC) \searrow M^{(n+1)}(d_{\a},\IC) \searrow
M^{(n+2)}(d_{\b},\IC) \searrow \cdots \searrow M^{(m)}(d_j,\IC)$.
\end{itemize}
~\\
The whole $\ca$ is an ideal which, by definition, is not primitive
since the trivial representation $\ca \rightarrow 0$ is not irreducible. 
Furthermore, the ideal $\{0\} \subset \ca$ is  primitive if and only
if $\ca$ has 
one irreducible faithful representation. This can be understood from the
Bratteli diagram in the following way. The set $\{0\}$ is not a subdiagram of
$\cd(\ca)$, being represented by the element $0$ of the matrix algebra of each
finite level, so that there is at least one element $a\in \ca$ not belonging to
the ideal $\{0\}$ at any level. Thus to check if $\{0\}$ is primitive, i.e. to
check property (iii) above, we have to see whether we can connect all the
points at a level $n$ to a {\it single} point at a level $m > n$. For example
this is the case for the diagram of Fig.~\ref{fi:veealg}. Later, we shall
construct examples in which $\{0\}$ is not a primitive ideal.
\medskip

Thus one can understand the topological properties of $Prim(\ca)$ at once 
from the
Bratteli diagram $\cd(\ca)$. This is  particularly simple if the algebra
admits only  a finite number of nonequivalent irreducible representations. In
this case $Prim(\ca)$ is a $T_0$-topological space with only a finite number of
points, hence a finite poset $P$. To reconstruct the latter we just
need to draw 
the Bratteli diagram $\cd(\ca)$ and find the subdiagrams that, according
to properties (i,ii,iii), correspond to primitive ideals. Then $P$
has as many points as the number of primitive ideals and the partial
order relation in $P$ is simply given by the inclusion relations that
exist among the primitive ideals.

\bexam
Again, consider the diagram of Fig.~\ref{fi:veealg}. The
corresponding AF algebra $\ca$ in (\ref{alvee}) contains only three nontrivial
ideals, whose diagrams are represented in
Fig.~\ref{fi:veealgid}(a,b,c). In this 
pictures the points belonging to the ideals are marked with a $``\clubsuit"$. 
It is not difficult to check that only $\ci_1$ and $\ci_2$ are
primitive ideals, 
since $\ci_3$ does not satisfy property (iii) above. Clearly, $\{0\}$ is
primitive and belongs to both $\ci_1$ and $\ci_2$ so that $Prim(\ca)$ is the
$\bigvee$ poset of Fig.~\ref{fi:bigvee}.
\begin{figure}[htb]
\begin{center}
\begin{picture}(360,160)(0,20)
\put(30,150){\circle*{4}}
\put(30,120){\circle*{4}}
\put(30,90){\circle*{4}}
\put(30,60){\circle*{4}}
\put(60,180){\circle*{4}}
\put(56,117){$\clubsuit$}
\put(56,87){$\clubsuit$}
\put(56,57){$\clubsuit$}
\put(86,147){$\clubsuit$}
\put(86,117){$\clubsuit$}
\put(86,87){$\clubsuit$}
\put(86,57){$\clubsuit$}
\put(30,150){\line(0,-1){30}}
\put(30,120){\line(0,-1){30}}
\put(30,90){\line(0,-1){30}}
\put(60,120){\line(0,-1){30}}
\put(60,90){\line(0,-1){30}}
\put(90,150){\line(0,-1){30}}
\put(90,120){\line(0,-1){30}}
\put(90,90){\line(0,-1){30}}
\put(30,60){\line(0,-1){10}}
\put(60,60){\line(0,-1){10}}
\put(90,60){\line(0,-1){10}}
\put(60,180){\line(1,-1){30}}
\put(30,150){\line(1,-1){30}}
\put(30,120){\line(1,-1){30}}
\put(30,90){\line(1,-1){30}}
\put(60,180){\line(-1,-1){30}}
\put(90,150){\line(-1,-1){30}}
\put(90,120){\line(-1,-1){30}}
\put(90,90){\line(-1,-1){30}}
\put(30,60){\line(1,-1){10}}
\put(90,60){\line(-1,-1){10}}
\put(146,147){$\clubsuit$}
\put(146,117){$\clubsuit$}
\put(146,87){$\clubsuit$}
\put(146,57){$\clubsuit$}
\put(180,180){\circle*{4}}
\put(176,117){$\clubsuit$}
\put(176,87){$\clubsuit$}
\put(176,57){$\clubsuit$}
\put(210,150){\circle*{4}}
\put(210,120){\circle*{4}}
\put(210,90){\circle*{4}}
\put(210,60){\circle*{4}}
\put(150,150){\line(0,-1){30}}
\put(150,120){\line(0,-1){30}}
\put(150,90){\line(0,-1){30}}
\put(180,120){\line(0,-1){30}}
\put(180,90){\line(0,-1){30}}
\put(210,150){\line(0,-1){30}}
\put(210,120){\line(0,-1){30}}
\put(210,90){\line(0,-1){30}}
\put(150,60){\line(0,-1){10}}
\put(180,60){\line(0,-1){10}}
\put(210,60){\line(0,-1){10}}
\put(180,180){\line(1,-1){30}}
\put(150,150){\line(1,-1){30}}
\put(150,120){\line(1,-1){30}}
\put(150,90){\line(1,-1){30}}
\put(180,180){\line(-1,-1){30}}
\put(210,150){\line(-1,-1){30}}
\put(210,120){\line(-1,-1){30}}
\put(210,90){\line(-1,-1){30}}
\put(150,60){\line(1,-1){10}}
\put(210,60){\line(-1,-1){10}}
\put(270,150){\circle*{4}}
\put(270,120){\circle*{4}}
\put(270,90){\circle*{4}}
\put(270,60){\circle*{4}}
\put(300,180){\circle*{4}}
\put(296,117){$\clubsuit$}
\put(296,87){$\clubsuit$}
\put(296,57){$\clubsuit$}
\put(330,150){\circle*{4}}
\put(330,120){\circle*{4}}
\put(330,90){\circle*{4}}
\put(330,60){\circle*{4}}
\put(270,150){\line(0,-1){30}}
\put(270,120){\line(0,-1){30}}
\put(270,90){\line(0,-1){30}}
\put(300,120){\line(0,-1){30}}
\put(300,90){\line(0,-1){30}}
\put(330,150){\line(0,-1){30}}
\put(330,120){\line(0,-1){30}}
\put(330,90){\line(0,-1){30}}
\put(270,60){\line(0,-1){10}}
\put(300,60){\line(0,-1){10}}
\put(330,60){\line(0,-1){10}}
\put(300,180){\line(1,-1){30}}
\put(270,150){\line(1,-1){30}}
\put(270,120){\line(1,-1){30}}
\put(270,90){\line(1,-1){30}}
\put(300,180){\line(-1,-1){30}}
\put(330,150){\line(-1,-1){30}}
\put(330,120){\line(-1,-1){30}}
\put(330,90){\line(-1,-1){30}}
\put(270,60){\line(1,-1){10}}
\put(330,60){\line(-1,-1){10}}
\put(52,15){$(a)$}
\put(172,15){$(b)$}
\put(292,15){$(c)$}
\put(75,170){$\ci_1$}
\put(195,170){$\ci_2$}
\put(315,170){$\ci_3$}
\put(45,30){$\vdots$}
\put(75,30){$\vdots$}
\put(165,30){$\vdots$}
\put(195,30){$\vdots$}
\put(285,30){$\vdots$}
\put(315,30){$\vdots$}
\end{picture}
\end{center}
\caption{\label{fi:veealgid} 
\protect{\footnotesize The three ideals of the algebra $\ca(\bigvee)$. }}
\end{figure}
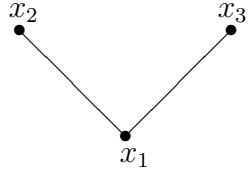
\begin{figure}[htb]
\begin{center}
\begin{picture}(320,80)(-70,40)
\put(80,40){\circle*{4}}
\put(40,80){\circle*{4}}
\put(120,80){\circle*{4}}
\put(80,40){\line(-1,1){40}}
\put(80,40){\line(1,1){40}}
\put(78,30){$x_1$}
\put(36,85){$x_2$}
\put(115,85){$x_3$}
\end{picture}
\end{center}
\caption{\label{fi:bigvee}
\protect{\footnotesize The poset $\bigvee$ as the primitive 
spectrum of the algebra $\ca(\bigvee)$. }}
\end{figure}\bigskip
\eexam

\subsxn{From Posets to Bratteli Diagrams}\label{se:bdp}

Under some rather mild hypotheses which are always verified in the cases
of posets, it is possible to reverse the construction of previous section and
thus construct an AF algebra starting from a poset.  
Such a reconstruction rests on another result by Bratteli \cite{Br2}, which
specifies the class of topological spaces which are the primitive ideal spaces
of AF algebras. Here, by using the  techniques of \cite{Br2} we shall
explain how 
to explicitly find an AF algebra $\ca$ (or rather its Bratteli diagram
$\cd(\ca)$) whose primitive ideal space is a given finite poset $P$. First we
will give the general construction and then discuss several examples.

Let $\{K_1,K_2,K_3,\ldots\}$ be the collection of all closed sets in the poset
$P$, with $K_1 = P$. To construct the $n$-th level of the Bratteli diagram
$\cd(\ca)$, we consider the subcollection $\ck_n = \{K_1,K_2,\ldots, K_n\}$ and
denote with $\ck_n'$ the smallest collection of (closed) sets in $P$
containing  
$\ck_n$ which is closed under union and intersection.  The collection $\ck_n$
determines a partition of the space $P$ by taking intersections and complements
of the sets $K_j \in \ck_n$ ($j=1,\ldots,n$). We denote with $Y(n,1),\; Y(n,2),
\; \ldots ,\;Y(n,k_n)$ the minimal sets of  such partition. Also, we write
$F(n,j)$ for the smallest set in the subcollection $\ck_n'$ which contains
$Y(n,j)$. Then, the diagram $\cd(\ca)$ is constructed as follows:
\begin{itemize}
\item[1.]
the $n$-th level of $\cd(\ca)$ has $k_n$ points, one for each set 
$Y(n,k), k=1, \cdots, k_n$;
\item[2.]
at the level $n$ of the diagram, the point which corresponds to $Y(n,i)$ is
linked to the point at the level $n+1$ corresponding to $Y(n+1,j)$ if and only
if $Y(n,i)\cap~F(n+1,j) \neq \emptyset$. In this case, the multiplicity of the
embedding is always 1.
\end{itemize}
\noindent
\bexam
To illustrate this construction, let us consider again the $\bigvee$
poset of Fig.~\ref{fi:bigvee}, $P=\{x_1, x_2, x_3 \}$. This topological
space contains four closed sets: 
\be
K_1=\{x_1, x_2, x_3 \}~, K_2 =\{x_2\}~, K_3 =\{x_3\}~, K_4 =\{x_2, x_3\} =
K_2
\cup K_3~.
\ee
Thus it is not difficult to check that:
\be
\begin{array}{llll}
\ck_1=\{K_1\} & \ck_1'=\{K_1\} & Y(1,1)=\{x_1, x_2, x_3\} & F(1,1)=K_1 \\ 
\ck_2=\{K_1,K_2\} & \ck_2'=\{K_1,K_2\} & Y(2,1)=\{x_2\}  & F(2,1)=K_2 \\
              &            & Y(2,2)=\{x_1, x_3\} & F(2,2)=K_1 \\
\ck_3=\{K_1,K_2,K_3\} & \ck_3'=\{K_1,K_2,K_3,K_4\}         
                                         & Y(3,1)=\{x_2\}  & F(3,1)=K_2 \\
                 &                       & Y(3,2)=\{x_1\}    & F(3,2)=K_1 \\
                 &                       & Y(3,3)=\{x_3\}   & F(3,3)=K_3 \\
\ck_4=\{K_1,K_2,K_3,K_4\} & \ck_4'=\{K_1,K_2,K_3,K_4\} 
                                     & Y(4,1)=\{x_2\}   & F(4,1)=K_2 \\
                  &                    & Y(4,2)=\{x_1\}    & F(4,2)=K_1 \\
                  &                    & Y(4,3)=\{x_3\}  & F(4,2)=K_3 \\
  &  \vdots & &
\end{array}
\ee
Since $P$ has only a finite number of points and hence a finite number
of closed sets, the partition of $P$ repeats itself after a certain
level ($n=3$ in this case). The corresponding diagram, obtained
through rules (1) 
and (2) above is readily found to coincide with the one in
Fig.~\ref{fi:veealg}. 
As we have said previously, such a diagram corresponds to the AF algebra 
$\ca(\bigvee) = \IC\cp_1 + \ck({\ch}) + \IC\cp_2$.
\eexam

It is a general fact that a Bratteli diagram describing any (finite) poset
`stabilizes' after a certain level $n_0$, namely it repeats itself . {}From
that level, both the number of points, as well as the embeddings from one level
to the next one, do not change. As we shall mention later, for the construction
of the $K$-theory groups, we only need the stable part of the Bratteli diagram
and we shall construct this part only in the remaining examples.

Now, the Bratteli diagram stabilizes at the level $n_0$ if the family
$\ck_{n_0}$ of closed sets we choose is such that it determines a
partition of the poset which distinguishes each point of the poset itself. In
particular, this is the case if we choose $n_0$ in such a manner that
$\ck_{n_0}$ contains all closed set. Then, each $Y(n_0,j)$ will contain a
single point of the poset and $F(n_0+1,j)$ will be the smallest closet set
containing $Y(n_0,j)$.

\bexam
The poset $\pinco$ of Fig.~\ref{fi:pinpos}.
Here $n_0=4$ and the stable partition is given by
\be
\begin{array}{ll}
Y(n_0,1)=\{x_3 \}  ~~&~~ F(n_0+1,1)=\{x_3 \} \\
Y(n_0,2)=\{x_1 \}  ~~&~~ F(n_0+1,2)=\{x_1, x_3, x_4 \} \\
Y(n_0,3)=\{x_2 \}  ~~&~~ F(n_0+1,3)=\{x_2, x_4 \} \\
Y(n_0,4)=\{x_4 \}  ~~&~~ F(n_0+1,4)=\{x_4 \} ~.
\end{array}
\ee
\begin{figure}[htb]
\begin{center}
\begin{picture}(230,60)(-110,-30)
\put(-30,30){\circle*{4}}
\put(30,30){\circle*{4}}
\put(-30,-30){\circle*{4}}
\put(30,-30){\circle*{4}}
\put(-30,30){\line(0,-1){60}}
\put(30,30){\line(0,-1){60}}
\put(30,30){\line(-1,-1){60}}
\put(-45,29){$x_3$}
\put(35,29){$x_4$}
\put(-45,-31){$x_1$}
\put(35,-31){$x_2$}
\end{picture}
\end{center}
\caption{\label{fi:pinpos} \protect{\footnotesize The poset $\pinco$.  }}
\end{figure}\bigskip

\noindent
The corresponding Bratteli diagram is shown in Fig.~\ref{fi:pinbra}.
The set $\{0\}$ is not an ideal. The algebra limit $\ca(\pinco)$ turns out to
be a subalgebra of bounded operators on the Hilbert space $\ch = \ch_1 \oplus
\ch_2 \oplus \ch_3$ with $\ch_i ,~ i = 1, 2, 3$ infinite dimensional Hilbert
spaces \cite{ELTfunctions},
\be
\ca(\pinco) = \IC\cp (\ch_1) \oplus \IC\cp (\ch_2 \oplus \ch_3)
             \oplus \ck(\ch_2 \oplus \ch_2) \oplus \ck(\ch_3) ~.
\label{pinalg}
\ee
Here, $\ck$ denotes the set of compact operators and $\cp$ orthogonal
projection. 
\begin{figure}[htb]
\begin{center}
\begin{picture}(120,120)(20,40)
\put(30,120){\circle*{4}}
\put(30,90){\circle*{4}}
\put(30,60){\circle*{4}}
\put(60,120){\circle*{4}}
\put(60,90){\circle*{4}}
\put(60,60){\circle*{4}}
\put(90,120){\circle*{4}}
\put(90,90){\circle*{4}}
\put(90,60){\circle*{4}}
\put(120,120){\circle*{4}}
\put(120,90){\circle*{4}}
\put(120,60){\circle*{4}}
\put(30,120){\line(0,-1){30}}
\put(30,90){\line(0,-1){30}}
\put(60,120){\line(0,-1){30}}
\put(60,90){\line(0,-1){30}}
\put(90,120){\line(0,-1){30}}
\put(90,90){\line(0,-1){30}}
\put(120,120){\line(0,-1){30}}
\put(120,90){\line(0,-1){30}}
\put(30,120){\line(1,-1){30}}
\put(30,90){\line(1,-1){30}}
\put(120,120){\line(-1,-1){30}}
\put(120,90){\line(-1,-1){30}}
\put(30,60){\line(0,-1){10}}
\put(30,60){\line(1,-1){10}}
\put(60,60){\line(0,-1){10}}
\put(90,60){\line(0,-1){10}}
\put(120,60){\line(-1,-1){10}}
\put(120,60){\line(0,-1){10}}
\put(120,60){\line(-2,-1){20}}
\put(30,120){\line(0,1){10}}
\put(60,120){\line(-1,1){10}}
\put(60,120){\line(0,1){10}}
\put(90,120){\line(0,1){10}}
\put(90,120){\line(1,1){10}}
\put(120,120){\line(0,1){10}}
\put(60,120){\line(2,1){20}}
\put(120,120){\line(-2,-1){60}}
\put(120,90){\line(-2,-1){60}}
\put(75,140){$\vdots$}
\put(75,30){$\vdots$}
\end{picture}
\end{center}
\caption{\label{fi:pinbra} 
\protect{\footnotesize The stable part of the Bratteli diagram of the poset
$\pinco$.}}
\end{figure}\bigskip
\eexam
\bexam
The poset $P_4(S^1)$ for the one-dimensional sphere in
Fig.~\ref{fi:cirhas}.  Here $n_0=4$ and the stable partition is given by
\be
\begin{array}{ll}
Y(n_0,1)=\{x_3 \}  ~~&~~ F(n_0+1,1)=\{x_3 \} \\
Y(n_0,2)=\{x_1 \}  ~~&~~ F(n_0+1,2)=\{x_1, x_3, x_4 \} \\
Y(n_0,3)=\{x_2 \}  ~~&~~ F(n_0+1,3)=\{x_2, x_3, x_4 \} \\
Y(n_0,4)=\{x_4 \}  ~~&~~ F(n_0+1,4)=\{x_4 \}~.
\end{array}
\ee
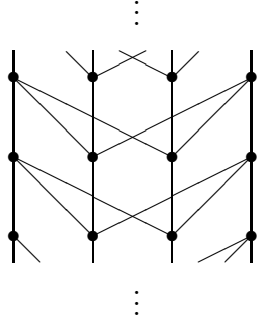
\begin{figure}[htb]
\begin{center}
\begin{picture}(120,120)(20,40)
\put(30,120){\circle*{4}}
\put(30,90){\circle*{4}}
\put(30,60){\circle*{4}}
\put(60,120){\circle*{4}}
\put(60,90){\circle*{4}}
\put(60,60){\circle*{4}}
\put(90,120){\circle*{4}}
\put(90,90){\circle*{4}}
\put(90,60){\circle*{4}}
\put(120,120){\circle*{4}}
\put(120,90){\circle*{4}}
\put(120,60){\circle*{4}}
\put(30,120){\line(0,-1){30}}
\put(30,90){\line(0,-1){30}}
\put(60,120){\line(0,-1){30}}
\put(60,90){\line(0,-1){30}}
\put(90,120){\line(0,-1){30}}
\put(90,90){\line(0,-1){30}}
\put(120,120){\line(0,-1){30}}
\put(120,90){\line(0,-1){30}}
\put(30,120){\line(1,-1){30}}
\put(30,90){\line(1,-1){30}}
\put(120,120){\line(-1,-1){30}}
\put(120,90){\line(-1,-1){30}}
\put(120,120){\line(-2,-1){60}}
\put(120,90){\line(-2,-1){60}}
\put(30,120){\line(2,-1){60}}
\put(30,90){\line(2,-1){60}}
\put(30,60){\line(0,-1){10}}
\put(30,60){\line(1,-1){10}}
\put(60,60){\line(0,-1){10}}
\put(90,60){\line(0,-1){10}}
\put(120,60){\line(-1,-1){10}}
\put(120,60){\line(0,-1){10}}
\put(120,60){\line(-2,-1){20}}
\put(30,120){\line(0,1){10}}
\put(60,120){\line(0,1){10}}
\put(60,120){\line(-1,1){10}}
\put(90,120){\line(0,1){10}}
\put(90,120){\line(1,1){10}}
\put(120,120){\line(0,1){10}}
\put(60,120){\line(2,1){20}}
\put(90,120){\line(-2,1){20}}
\put(75,140){$\vdots$}
\put(75,30){$\vdots$}
\end{picture}
\end{center}
\caption{\label{fi:cirbra} 
\protect{\footnotesize The stable part of the Bratteli diagram
for the circle poset $P_4(S^1)$  . }}
\end{figure}\bigskip

\noindent
The corresponding Bratteli diagram is shown in Fig.~\ref{fi:cirbra}.
The set $\{0\}$ is not an ideal. The algebra limit $\ca(P_4(S^1))$ turns out to
be a subalgebra of bounded operators on the Hilbert space $\ch = \ch_1 \oplus
\cdots \oplus \ch_4$~, with $\ch_i ,~ i = 1, \dots, 4$ infinite
dimensional Hilbert spaces \cite{ELTfunctions},
\be
\ca(P_4(S^1)) = \IC\cp (\ch_1 \oplus \ch_4) \oplus \IC\cp (\ch_2 \oplus \ch_3)
\oplus \ck (\ch_1 \oplus \ch_2) \oplus \ck (\ch_3 \oplus \ch_4) ~.
\label{ciralg}
\ee
Here $\ck$ denotes the set of compact operators and $\cp$ orthogonal
projection. 
\eexam
\bexam
The poset $P_6(S^2)$ for the two-dimensional sphere in Fig.~\ref{fi:sphpos}. 
Here $n_0=6$ and the stable partition is given by
\be
\begin{array}{ll}
Y(n_0,1)=\{x_5 \}  ~~&~~ F(n_0+1,1)=\{x_5 \}\\
Y(n_0,2)=\{x_3 \}  ~~&~~ F(n_0+1,2)=\{x_3, x_5, x_6 \} \\
Y(n_0,3)=\{x_1 \}  ~~&~~ F(n_0+1,3)=\{x_1, x_3, x_4, x_5, x_6, \} \\
Y(n_0,4)=\{x_2 \}  ~~&~~ F(n_0+1,4)=\{x_2, x_3, x_4, x_5, x_6 \} \\
Y(n_0,5)=\{x_4 \}  ~~&~~ F(n_0+1,5)=\{x_4, x_5, x_6 \} \\
Y(n_0,6)=\{x_6 \}  ~~&~~ F(n_0+1,6)=\{x_6 \} ~.
\end{array}
\ee
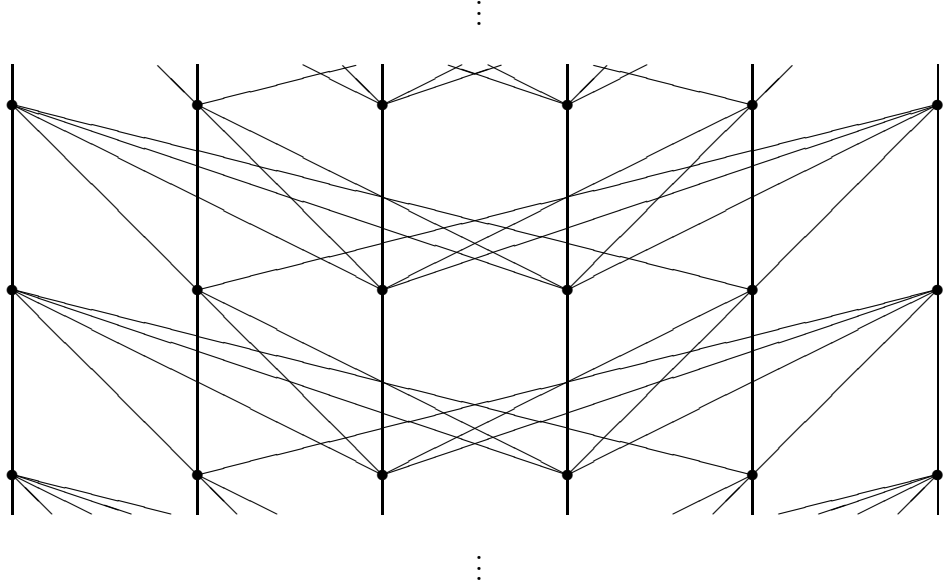
\begin{figure}[htb]
\begin{center}
\begin{picture}(600,220)(-230,60)
\put(-175,100){\circle*{4}}
\put(-175,170){\circle*{4}}
\put(-175,240){\circle*{4}}
\put(-105,100){\circle*{4}}
\put(-105,170){\circle*{4}}
\put(-105,240){\circle*{4}}
\put(-35,100){\circle*{4}}
\put(-35,170){\circle*{4}}
\put(-35,240){\circle*{4}}
\put(175,100){\circle*{4}}
\put(175,170){\circle*{4}}
\put(175,240){\circle*{4}}
\put(105,100){\circle*{4}}
\put(105,170){\circle*{4}}
\put(105,240){\circle*{4}}
\put(35,100){\circle*{4}}
\put(35,170){\circle*{4}}
\put(35,240){\circle*{4}}
\put(-175,170){\line(0,-1){70}}
\put(-175,240){\line(0,-1){70}}
\put(-105,170){\line(0,-1){70}}
\put(-105,240){\line(0,-1){70}}
\put(-35,170){\line(0,-1){70}}
\put(-35,240){\line(0,-1){70}}
\put(175,170){\line(0,-1){70}}
\put(175,240){\line(0,-1){70}}
\put(105,170){\line(0,-1){70}}
\put(105,240){\line(0,-1){70}}
\put(35,170){\line(0,-1){70}}
\put(35,240){\line(0,-1){70}}
\put(-175,100){\line(0,-1){15}}
\put(-105,100){\line(0,-1){15}}
\put(-35,100){\line(0,-1){15}}
\put(175,100){\line(0,-1){15}}
\put(105,100){\line(0,-1){15}}
\put(35,100){\line(0,-1){15}}
\put(-175,240){\line(0,1){15}}
\put(-105,240){\line(0,1){15}}
\put(-35,240){\line(0,1){15}}
\put(175,240){\line(0,1){15}}
\put(105,240){\line(0,1){15}}
\put(35,240){\line(0,1){15}}
\put(-175,170){\line(1,-1){70}}
\put(-175,170){\line(2,-1){140}}
\put(-175,170){\line(3,-1){210}}
\put(-175,170){\line(4,-1){280}}
\put(-175,240){\line(1,-1){70}}
\put(-175,240){\line(2,-1){140}}
\put(-175,240){\line(3,-1){210}}
\put(-175,240){\line(4,-1){280}}
\put(175,170){\line(-1,-1){70}}
\put(175,170){\line(-2,-1){140}}
\put(175,170){\line(-3,-1){210}}
\put(175,170){\line(-4,-1){280}}
\put(175,240){\line(-1,-1){70}}
\put(175,240){\line(-2,-1){140}}
\put(175,240){\line(-3,-1){210}}
\put(175,240){\line(-4,-1){280}}
\put(-105,170){\line(1,-1){70}}
\put(-105,170){\line(2,-1){140}}
\put(-105,240){\line(1,-1){70}}
\put(-105,240){\line(2,-1){140}}
\put(105,170){\line(-1,-1){70}}
\put(105,170){\line(-2,-1){140}}
\put(105,240){\line(-1,-1){70}}
\put(105,240){\line(-2,-1){140}}
\put(-175,100){\line(1,-1){15}}
\put(-175,100){\line(2,-1){30}}
\put(-175,100){\line(3,-1){45}}
\put(-175,100){\line(4,-1){60}}
\put(175,100){\line(-1,-1){15}}
\put(175,100){\line(-2,-1){30}}
\put(175,100){\line(-3,-1){45}}
\put(175,100){\line(-4,-1){60}}
\put(-105,100){\line(1,-1){15}}
\put(-105,100){\line(2,-1){30}}
\put(105,100){\line(-1,-1){15}}
\put(105,100){\line(-2,-1){30}}
\put(-35,240){\line(-1,1){15}}
\put(-35,240){\line(-2,1){30}}
\put(-35,240){\line(2,1){30}}
\put(-35,240){\line(3,1){45}}
\put(35,240){\line(1,1){15}}
\put(35,240){\line(2,1){30}}
\put(35,240){\line(-3,1){45}}
\put(35,240){\line(-2,1){30}}
\put(-105,240){\line(-1,1){15}}
\put(-105,240){\line(4,1){60}}
\put(105,240){\line(1,1){15}}
\put(105,240){\line(-4,1){60}}
\put(0,60){$\vdots$}
\put(0,270){$\vdots$}
\end{picture}
\end{center}
\caption{\label{fi:sphbra} 
\protect{\footnotesize The stable part of the Bratteli diagram
for the sphere poset $P_6(S^2)$.  }}
\end{figure}\bigskip

\noindent
The corresponding Bratteli diagram is shown in Fig.~\ref{fi:sphbra}.
The set $\{0\}$ is not an ideal. The algebra limit $\ca(P_6(S^2))$ turns out to
be a subalgebra of bounded operators on the Hilbert space $\ch = \ch_1 \oplus
\cdots \oplus \ch_8$ with $\ch_i ,~ i = 1, \dots, 8$ infinite dimensional
Hilbert spaces \cite{ELTfunctions},
\bea
&&\ca(P_6(S^2)) = \IC\cp \left(\ch_5 \otimes (\ch_1 \oplus \ch_4) \oplus
                          \ch_8 \otimes (\ch_2 \oplus \ch_3) \right)
\nonumber \\
&& ~~~~~~~~~~~~~~ \oplus \IC\cp \left(\ch_6 \otimes (\ch_1 \oplus \ch_4) \oplus
                          \ch_7 \otimes (\ch_2 \oplus \ch_3) \right)
\nonumber \\
&& \oplus
\left[\ck (\ch_5 \oplus \ch_6) \otimes \IC\cp (\ch_1 \oplus \ch_4) \right]
\oplus
\left[\ck (\ch_7 \oplus \ch_8) \otimes \IC\cp (\ch_2 \oplus \ch_3) \right]
\nonumber \\
&& \oplus
\ck \left[ \ch_1 \otimes (\ch_5 \oplus \ch_6) \oplus
                         \ch_2 \otimes (\ch_7 \oplus \ch_8) \right]
\oplus
\ck \left[ \ch_3 \otimes (\ch_7 \oplus \ch_8) \oplus
                         \ch_4 \otimes (\ch_5 \oplus \ch_6) \right]
 \nonumber \\
 ~ \label{sphalg}
\eea
Here $\ck$ denotes the set of compact operators and $\cp$ orthogonal
projection. 
\eexam

\sxn{Projective modules of finite type and $K$-theory}\label{se:pbk}

Given an algebra $\ca$ playing the role of the algebra of continuous functions
on some noncommutative space, the analogue of vector bundles is provided by the
notion of {\it projective module of finite type} (or {\it finite projective
module}) over $\ca$. Indeed, by the Serre-Swan theorem \cite{SW}, locally
trivial, finite-dimensional complex vector bundles over a compact Hausdorff
space $M$ are in one to one correspondence with finite projective modules over
the algebra $\ca = \cc(M)$. To the vector bundle $E$ one associates the
$\cc(M)$-module $\ce=\G(M,E)$ of continuous sections of
$E$. Conversely, if $\ce$ 
is a finite projective modules over $\cc(M)$, the fiber $E_m$ of the associated
bundle $E$ over the point
$m \in M$ is
\be
E_m = \ce / \ce \ci_m~,
\ee
where the ideal $\ci_m \subset \cc(M)$, corresponding to the point $m \in M$,
is given by \cite{Co,VG}
\be\ci_m = ker\{ \c_m : \cc(M)
\ra \IC~; ~ \c_m(f) =  f(m) \} = \{ f \in \cc(M) ~|~ f(m) = 0 \}~.
\ee

As we shall see, isomorphism and stable isomorphism have a meaning in
the context of finite projective modules over a \cstar $\ca$ and the
group $K_0(\ca)$ will be the group of (stable) isomorphism classes of
finite projective (right) modules over $\ca$.

Given any finite projective right module $\ce$ over $\ca$, there exists 
an integer $N$ together with an idempotent $p \in \IM(N, \ca)$ ($N \times N$
matrix with entries in $\ca$ and $p^2=p$), and an isomorphism   of
$\ce$ with the 
right $\ca$-module \cite{Co}
\be
p\ca^N = \{ \x = (\x_1, \dots, \x_N)~; ~\x_i \in \ca~, ~p\x =\x \}~.
\ee
In fact, we shall restrict to the class of {\it hermitian} modules
which correspond to projectors, namely idempotents $p$ obeying the additional
condition $p^* = p$, the operation $^*$ being the composition of the $^*$
operation in the algebra $\ca$ with usual matrix transposition.

We shall now give few fundamentals of the $K$-theory of \cstars having in mind
mainly AF algebras \cite{W-O}. Two projectors $p, q \in \IM(N, \ca)$ are
equivalent if there exists a matrix $u \in \IM(N, \ca)$ such that
$p=u^* u$ and $q=u u^*$. In order to be able to add equivalence classes of
projectors, one considers all finite matrix algebras over $\ca$ at the same
time, by considering $\IM(\infty, \ca)$ which is the non complete $^*$-algebra
obtained as inductive limit of finite matrices
\fn{The completion of $\IM(\infty, \ca)$ is $\ca \otimes \ck$, with $\ck$ the
algebra of compact operators on $l_2$. The algebra $\ca \otimes \ck$ is also
called the stabilization of $\ca$.},
\bea
&& \IM(\infty, \ca) = \bigcup_{n=1}^{\infty} \IM(n, \ca)~, \nonumber \\
&& \f : \IM(n, \ca) \ra \IM(n+1, \ca)~, ~a \mapsto
\f(a) =
\left[
\begin{array}{cc}
a & 0 \\
0 & 0
\end{array}
\right]~.
\eea

Now, two projectors $p, q \in \IM(\infty, \ca)$ are said to be equivalent, $p
\sim q$, when there exists an $u \in \IM(\infty, \ca)$ such that $p=u^* u$ and
$q=u u^*$. The set $V(\ca)$ of equivalence classes $[ ~\cdot~]$ is made into an
abelian semigroup by defining an {\it addition} by
\be
[p] + [q] =: [
\left[
\begin{array}{cc}
p & 0 \\
0 & q
\end{array}
\right]
] ~.
\ee

The group $K_0(\ca)$ is the universal canonical group associated with the
semigroup $V(\ca)$ and may be defined as
\bea
&& K_0(\ca) =: V(\ca) \times V(\ca) / \sim ~, ~~([p], [q]) \sim ([p'], [q'])~,
\nonumber \\
&& \iff~~ ~{\rm there~exists}~ [r] \in V(\ca)
~~{\rm s.t.}~~ [p] + [q'] + [r] = [p'] + [q] + [r]~. \label{keya}
\eea
The extra $[r]$ is necessary to get transitivity and make $\sim$ an equivalence
relation. This is the reason why one is classifying only stable classes.
{}From definition (\ref{keya}), an equivalence class $[([p], [q])] \in
K_0(\ca)$ 
can also be written as a formal difference $[p] - [q]$.

There is a natural homomorphism
\be
\k_{\ca} : V(\ca) \ra K_0(\ca)~,~~\k_{\ca}([p]) =: [p] - [0]  \label{homkey}
\ee
This map is injective if and only if the addition in $V(\ca)$ has
cancellations, namely if and only if $[p] + [r] = [q] + [r]~ \Rightarrow
[p] = [q]$.

While for a generic $\ca$, the semigroup $V(\ca)$ has no
cancellations, for AF algebras this happens to be the case. By defining
\be
K_{0+}(\ca) =: \k_{\ca}(V(\ca))~, \label{keya+}
\ee
the couple $(K_0(\ca), K_{0+}(\ca))$ becomes, for an AF algebra $\ca$,
an {\it ordered group} with $K_{0+}(\ca)$ the {\it positive cone}, namely
one has that 
\bea
&& K_{0+}(\ca) \ni 0~, \nonumber \\
&& K_{0+}(\ca) - K_{0+}(\ca) = K_0(\ca)~, \nonumber \\
&& K_{0+}(\ca) \cap (- K_{0+}(\ca)) = 0~.
\eea
For a generic (unital) algebra the last property is not true and the couple
$(K_0(\ca), K_{0+}(\ca))$ is not an ordered group.

There is another $K$-group, $K_1$, which is constructed from unitaries or 
invertibles in $\IM(\infty, \ca)$. It turns out, however, that such a group is
trivial for AF  algebras, namely $K_1(\ca)=0$ for any AF algebras
\cite{W-O}. We 
shall  not mention it anymore in the rest of the paper.

The construction of the $K$-theory of AF algebras is made easy by the following
results whose proofs are given in \cite{W-O}.

\bprop\label{pr:ind}
If $\a : \ca \ra \cb$ is a homomorphism of \cstars, then the induced map
\be
\a_* : V(\ca) \ra V(\cb)~, ~~\a_*([a_{ij}]) =: [\a(a_{ij})]~,
\label{induc}
\ee
is a well defined homomorphism of semigroups. Moreover, from universality,
$\a_*$ extends to a group homomorphism
\be
\a_* = K_0(\ca) \ra K_0(\cb)~. \label{keyhom}
\ee
\eprop
\bprop\label{pr:key}
If $\ca$ is the inductive limit of a directed system $\{\ca_i, \F_{ij}\}_I$ of
\cstars, then $\{K_0(\ca_{i}), \F_{ij*}\}_I$ is a directed system of groups and
one can exchange the limits,
\be
K_0(\ca) = K_0(\lim_{\ra} \ca_i) = \lim_{\ra} K_0(\ca_i)~.
\ee
Moreover, if $\ca$ is an AF algebra, then $K_0(\ca)$ is an ordered group with
positive cone given by the limit of a directed system of semigroups
\be
K_{0+}(\ca) = K_{0+}(\lim_{\ra} \ca_i) = \lim_{\ra} K_{0+}(\ca_i)~.
\ee
\eprop
\bprop\label{pr:keyaf}
With  $k_{\ca}, k_{\cb}$ integer numbers, let $\ca$ and $\cb$ be the sum of
$k_{\ca}$ and $k_{\cb}$ matrix algebras respectively, 
\bea
&&\ca = \IM(p_1, \IC) \oplus \IM(p_2, \IC) \oplus \cdots
   \oplus \IM(p_{k_{\ca}}, \IC)~, \nonumber \\
&&\cb = \IM(q_1, \IC) \oplus \IM(q_2, \IC) \oplus \cdots
   \oplus \IM(q_{k_{\cb}}, \IC)~. \nonumber \\
\eea
Then, any homomorphism $\a : \ca \ra \cb$ can be written as the direct sum of
the representations $\a_j : \ca \ra \IM(q_j, \IC) \simeq \cb(\IC^{q_j})~,
j=1,\cdots,k_{\cb}$. If $\p_{ij}$ is the unique irreducible representation of 
$\IM(p_i, \IC)$ in $\cb(\IC^{q_j})$, then $\a_j$ breaks into a direct sum of
the $\p_{ij}$. Furthermore, let $m_{ij}$ be the non-negative integer denoting
the multiplicity of $\p_{ij}$ in this sum. Then the induced homomorphism,
$\a_* = K_0(\ca) \ra K_0(\cb)$, is given by the $q_{k_{\ca}} \times
p_{k_{\cb}}$ 
matrix $(m_{ij})$. 
\eprop

We end this section by mentioning that $K$-theory has been proved
\cite{El} to be a complete invariant which distinguish among AF algebras if one
add to the ordered group $(K_0(\ca), K_{0+}(\ca))$ the notion of {\it scale},
the latter being defined for any \cstar $\ca$ as
\be
\S\ca =: \{[p]~, p ~~{\rm a~ projector~ in~ \ca } \}~.
\ee
AF algebras are completely determined, up to isomorphism, by their {\it scaled
ordered} groups, namely by triple $(K_0, K_{0+}, \S)$. The key is the fact that
scale preserving isomorphisms between the ordered groups $(K_0, K_{0+}, \S)$ of
two AF algebras are nothing but $K$-theoretically induced maps
(\ref{keyhom}) of 
isomorphisms between the AF algebras themselves.

\subsxn{The $K$-theory of noncommutative lattices}

The starting point to compute the ordered group $(K_0, K_{0+})$ for a
poset is the fact that, for an AF algebra given as in \ref{af}, the
group $(K_0(\ca), K_{0+}(\ca))$ is obtained by Proposition \ref{pr:key} as the
inductive limit of the sequence of groups/semigroups
\bea
&&
K_0(\ca_1) \hookrightarrow K_0(\ca_2) \hookrightarrow K_0(\ca_3)
\hookrightarrow 
\cdots
\label{kappa}\\
&&
K_{0+}(\ca_1) \hookrightarrow K_{0+}(\ca_2) \hookrightarrow K_{0+}(\ca_3)
\hookrightarrow
\cdots
\label{kappa+}
\eea
The inclusions
\be
T_n : K_0(\ca_n) \hookrightarrow K_0(\ca_{n+1})~,~~
T_n : K_{0+}(\ca_n) \hookrightarrow K_{0+}(\ca_{n+1})~, \label{emb}
\ee
are easily obtained from the inclusions
$\ca_n ~{\buildrel I_n \over \hookrightarrow}~ \ca_{n+1}$, being indeed the
corresponding induced maps as in (\ref{keyhom}). As sets we have that
\bea
K_0(\ca) = \{ (k_n)_{n \in \IN}~, k_n \in K_0(\ca_n) ~|~ \exists N_0~:~
k_{n+1} = T_n(k_n)~,~ n>N_0 \}~, \label{key} \\
K_{0+}(\ca) = \{ (k_n)_{n \in \IN}~, k_n \in K_{0+}(\ca_n) ~|~
\exists N_0~:~ k_{n+1} = T_n(k_n)~,~ n>N_0 \}~. \label{key+}
\eea
The structure of (abelian) group/semigroup is inherited pointwise from the
addition in the groups/semigroups in the sequences
(\ref{kappa}), (\ref{kappa+}).

Furthermore, for any $d$, the algebra of matrices $\IM(d, \IC)$ has $K$-theory
given by $(K_0, K_{0+}) = (\IZ, \IZ_+)$, $\IZ$ being the group of integer
numbers and $\IZ_+$ the semigroups of natural numbers (including
$0$). Hence,
all terms in the sequences (\ref{kappa}), (\ref{kappa+}), are direct sums of
copies of $\IZ$ or $\IZ_+$.

As mentioned in Section \ref{se:bdp}, the Bratteli diagrams that describe
(finite) posets have the property that starting from a certain level $n_0$
(which is less or equal than the number of closed sets in the poset),
the number 
of points in any diagram, as well as the embeddings from one level to the next
one, does not change. This simplifies the calculation of $(K_0, K_{0+})$ 
because, 
\bea
&&K_0(\ca_{n_0}) = K_0(\ca_{n_0 + 1}) = K_0(\ca_{n_0 + 2}) = \cdots =
\IZ^{\oplus k_{n_0}}~, \\
&&K_{0+}(\ca_{n_0}) = K_{0+}(\ca_{n_0 + 1}) = K_{0+}(\ca_{n_0 + 2})
= \cdots = \IZ_+^{\oplus k_{n_0}}~,
\eea
where $k_{n_0}$ is the number of points in the Bratteli diagram from the level
$n_0$ on. Furthermore, the integer valued matrices $T_n$  in
(\ref{emb}) are all 
equal for $n>n_0$.
To find the group $(K_0(\ca), K_{0+}(\ca))$ one has just to study  
the limit for $n \ra \infty$ of the inclusions
\bea
&& T_n : \IZ^{\oplus k_{n_0}} \hookrightarrow \IZ^{\oplus k_{n_0}} ~,
\label{keyinc}\\
&& T_n : \IZ_+^{\oplus k_{n_0}} \hookrightarrow \IZ_+^{\oplus k_{n_0}} ~.
\label{keyinc+}
\eea
We infer from Prop.~\ref{pr:keyaf} that for AF algebras the maps
(\ref{keyinc}), 
(\ref{keyinc+}) are always inclusions.
In fact, for (finite) posets, the map (\ref{keyinc}) is always a
bijection. As a 
consequence, for a poset $P$ with $k_{n_0}$ points in the stable part of the
corresponding Bratteli diagram, and associated algebra $\ca_{k_{n_0}}(P)$, we
shall have that
\be
K_0(P) = \IZ^{\oplus k_{n_0}}~.
\ee
The map (\ref{keyinc+}) will not be in general a bijection.

\bigskip

\begin{figure}[htb]
\begin{center}
\begin{picture}(120,180)(20,50)
\put(60,210){\circle*{4}}
\put(30,180){\circle*{4}}
\put(30,150){\circle*{4}}
\put(30,120){\circle*{4}}
\put(30,90){\circle*{4}}
\put(30,60){\circle*{4}}
\put(90,180){\circle*{4}}
\put(90,150){\circle*{4}}
\put(90,120){\circle*{4}}
\put(90,90){\circle*{4}}
\put(90,60){\circle*{4}}
\put(60,210){\line(-1,-1){30}}
\put(60,210){\line(1,-1){30}}
\put(30,180){\line(0,-1){30}}
\put(30,150){\line(0,-1){30}}
\put(30,120){\line(0,-1){30}}
\put(30,90){\line(0,-1){30}}
\put(30,60){\line(0,-1){10}}
\put(30,180){\line(2,-1){60}}
\put(30,150){\line(2,-1){60}}
\put(30,120){\line(2,-1){60}}
\put(30,90){\line(2,-1){60}}
\put(30,150){\line(2,1){60}}
\put(30,120){\line(2,1){60}}
\put(30,90){\line(2,1){60}}
\put(30,60){\line(2,1){60}}
\put(58,214){1}
\put(21,178){1}
\put(21,148){2}
\put(21,118){3}
\put(21,88){5}
\put(21,58){8}
\put(94,178){1}
\put(94,148){1}
\put(94,118){2}
\put(94,88){3}
\put(94,58){5}
\put(60,40){$\vdots$}
\end{picture}
\end{center}
\caption{\label{fi:penbra} 
\protect{\footnotesize The Bratteli diagram for the algebra of
the Penrose tiling.   }}
\end{figure}
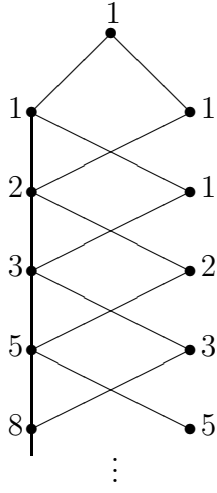\bigskip

We shall illustrate the construction of the 
$K$-groups with the Penrose Tiling AF
algebra. Although this algebra is quite far from being postliminal, since there
are infinite non equivalent representations with the same kernel (the only
primitive ideal), the construction of its $K$-theory is illuminating. The
corresponding Bratteli diagram is shown in Fig.~\ref{fi:penbra}.
At each level, the algebra is given by \cite{Co}
\be
\ca_n = \IM(d_n, \IC) \oplus \IM(d'_n, \IC)~, ~~ n \geq 1~, \label{pt}
\ee
with inclusion $\ca_n \hookrightarrow \ca_{n+1}$,
\be\label{ptinc}
\left[
\begin{array}{cc}
A & 0 \\
0 & B
\end{array}
\right] \ra
\left[
\begin{array}{ccc}
A & 0 & 0 \\
0 & B & 0 \\
0 & 0 & A
\end{array}
\right]~,
~~A \in \IM(d_n, \IC)~, ~~B \in \IM(d'_n, \IC)~.
\ee
After the second level we have then
\be
K_0(\ca_n) = \IZ \oplus \IZ~, ~~~K_{0+}(\ca_n) = \IZ_+ \oplus \IZ_+~.
\ee
The inclusion
(\ref{ptinc}) gives for the dimensions
\bea
&& d_{n+1} = d_n + d'_n~, \nonumber \\
&& d'_{n+1} = d_n~. \label{ptdim}
\eea
while the inclusions (\ref{keyinc}), (\ref{keyinc+}) are both
represented by the 
integer valued matrix
\be
T =
\left[
\begin{array}{cc}
1 & 1 \\
1 & 0
\end{array}
\right]~. \label{ptmat}
\ee
The action of the matrix (\ref{ptmat}) can be represented pictorially as
in Fig.~\ref{fi:ptmatpic} where the couple $(a, b)$ $(a', b')$ are both in 
$\IZ \oplus \IZ$ or $\IZ_+ \oplus \IZ_+$.
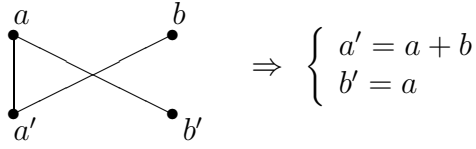
\begin{figure}[htb]
\begin{center}
\begin{picture}(200,40)(0,60)
\put(30,90){\circle*{4}}
\put(30,60){\circle*{4}}
\put(90,90){\circle*{4}}
\put(90,60){\circle*{4}}
\put(30,90){\line(0,-1){30}}
\put(30,90){\line(2,-1){60}}
\put(30,60){\line(2,1){60}}
\put(30,94){$a$}
\put(30,50){$a'$}
\put(90,94){$b$}
\put(94,50){$b'$}
\put(120,75){$\Rightarrow$}
\put(140,75)
{$
\left\{
\begin{array}{l}
a'=a+b \\
b'=a
\end{array}
\right.
$}
\end{picture}
\end{center}
\caption{\label{fi:ptmatpic} 
\protect{\footnotesize The action of the inclusion $T$. }}
\end{figure}\bigskip

Finally, we can construct the $K$-theory groups.
\begin{itemize}
\item[1)] The group $K_0(\ca)$ is given by
\be
K_0(\ca) = \IZ \oplus \IZ~. \label{ptkey}
\ee
This follows immediately from the fact that the matrix $T$ in (\ref{ptmat}) is
invertible over the integer, its inverse being
\be
T^{-1} =
\left[
\begin{array}{cc}
0 & 1 \\
1 & -1
\end{array}
\right]~. \label{ptmatinv}
\ee
Now, from the definition of inductive limit we have that,
\be
K_0(\ca) = \{ (k_n)_{n \in \IN}~, k_n \in K_0(\ca_n) ~|~ \exists N_0 ~:~
k_{n+1} = T (k_n)~, ~n>N_0 \}.
\ee
In addition, $T$ being a bijection, for any $k_{n+1} \in K_0(\ca_{n+1})$, there
exist an unique $k_n \in K_0(\ca_n)$ such that $k_{n+1} = T k_n$.  Thus,
$K_0(\ca) = K_0(\ca_n) = \IZ \oplus \IZ$.\\

\item[2)] The group $K_{0+}(\ca)$ is given by
\be
K_{0+}(\ca) = 
\{ (a, b) \in \IZ \oplus \IZ~ : {1+\sqrt{5} \over 2} a + b \geq 0 \}~.
\label{ptkey+}
\ee
Since $T$ is not invertible over $\IZ_+$, $K_{0+}(\ca) \not= \IZ_+
\oplus \IZ_+$. To construct $K_{0+}(\ca)$, we study the image
$T(K_{0+}(\ca_n))$ 
in $K_{0+}(\ca_{n+1})$. It is easily found to be 
\bea
T(K_{0+}(\ca_n)) & = & \{(a_{n+1}, b_{n+1}) \in \IZ_+ \oplus \IZ_+~ : 
a_{n+1} \geq b_{n+1} \} \nonumber \\
                 &\not=& K_{0+}(\ca_{n+1}) ~.
\eea
Now, $T$ being injective, $T(K_{0+}(\ca_n)) = T(\IZ_+ \oplus \IZ_+)
\simeq \IZ_+ \oplus \IZ_+$. The inclusion of $T(K_{0+}(\ca_n))$ into
$K_{0+}(\ca_{n+1})$ is shown in  Fig.~\ref{fi:penrosekeyplus}. 
If we identify the subset $T(K_{0+}(\ca_n)) \subset K_{0+}(\ca_{n+1})$ with
$K_{0+}(\ca_n)$,  we can think of $T^{-1}(K_{0+}(\ca_{n+1}))$ as a subset
of $\IZ \oplus \IZ$ and of $T^{-1}(K_{0+}(\ca_n))$ as the standard
positive cone 
$\IZ_{+} \oplus \IZ_{+}$. The result is shown in
Fig.~\ref{fi:penrosekeyminus}.  
Next iteration, namely $T^{-2}(K_{0+}(\ca_n))$ is shown in
Fig.~\ref{fi:penrosekeyminus2}. \\
{}From definition (\ref{key+}), by going to the
limit we shall have $K_{0+}(\ca) = \lim_{m \ra \infty} T^{-m}(\IZ_+ \oplus
\IZ_+)$ and the limit will be a subset of $ \IZ \oplus \IZ$ since $T$ is
invertible only over $\IZ$. The limit can be easily found 
\fn{We owe the following method to J. Varilly.}.
{}From the defining relation $F_{m+1} = F_m + F_{m-1}, m\geq1, $ for the 
Fibonacci numbers  (with $F_0=0, F_1=1$), it follows that
\be
T^{-m} =
(-1)^{m} \left[
\begin{array}{cc}
F_{m-1} & -F_m \\
-F_m & F_{m+1}
\end{array}
\right]~. \label{ptmatinvm}
\ee
Therefore, $T^{-m}$ takes the positive axis $\{(a,0) : a \geq 0\}$ to an
half-line of slope $-F_m / F_{m-1}$, and the positive axis $\{(0,b) : b \geq
0\}$  to an half-line of slope $-F_{m+1} / F_m$. Thus the positive cone
$\IZ_+ \oplus \IZ_+$ opens into a fan-shaped wedge bordered by these two
half-lines. Any integer coordinate point within the wedge comes from an
integer coordinate point in the original positive cone. 
Since $\lim_{m \ra \infty} F_{m+1} / F_m = {1+\sqrt{5} \over 2}$, the limit
cone is just the half-space $\{ (a, b) \in \IZ \oplus \IZ~ :
{1+\sqrt{5} \over 2} a + b \geq 0 \}$~. Every integer coordinate point in it
belongs to some intermediate wedge and so lies in $K_{0+}(\ca)$. The latter
is shown in  Fig.~\ref{fi:penkeyzerplu}.
\end{itemize}

\begin{figure}[htb]
\begin{center}
\begin{picture}(600,130)(-225,40)
\put(-180,50){\vector(1,0){130}}
\put(-170,40){\vector(0,1){130}}
\put(30,50){\vector(1,0){130}}
\put(40,40){\vector(0,1){130}}
\put(40,50){\vector(1,1){110}}
\put(-70,40){$a_n$}
\put(-183,150){$b_n$}
\put(140,40){$a_{n+1} = b_{n}$}
\put(18,150){$b_{n+1}$}
\put(152,151){$a_n$}
\put(-175,47){$\times$}
\put(-175,67){$\times$}
\put(-175,87){$\times$}
\put(-175,107){$\times$}
\put(-175,127){$\times$}
\put(-155,47){$\times$}
\put(-155,67){$\times$}
\put(-155,87){$\times$}
\put(-155,107){$\times$}
\put(-155,127){$\times$}
\put(-135,47){$\times$}
\put(-135,67){$\times$}
\put(-135,87){$\times$}
\put(-135,107){$\times$}
\put(-135,127){$\times$}
\put(-115,47){$\times$}
\put(-115,67){$\times$}
\put(-115,87){$\times$}
\put(-115,107){$\times$}
\put(-115,127){$\times$}
\put(-95,47){$\times$}
\put(-95,67){$\times$}
\put(-95,87){$\times$}
\put(-95,107){$\times$}
\put(-95,127){$\times$}
\put(-20,100){$\stackrel{T}{\longrightarrow}$}
\put(40,50){\circle*{3}}
\put(40,70){\circle*{3}}
\put(40,90){\circle*{3}}
\put(40,110){\circle*{3}}
\put(40,130){\circle*{3}}
\put(60,50){\circle*{3}}
\put(60,70){\circle*{3}}
\put(60,90){\circle*{3}}
\put(60,110){\circle*{3}}
\put(60,130){\circle*{3}}
\put(80,50){\circle*{3}}
\put(80,70){\circle*{3}}
\put(80,90){\circle*{3}}
\put(80,110){\circle*{3}}
\put(80,130){\circle*{3}}
\put(100,50){\circle*{3}}
\put(100,70){\circle*{3}}
\put(100,90){\circle*{3}}
\put(100,110){\circle*{3}}
\put(100,130){\circle*{3}}
\put(120,50){\circle*{3}}
\put(120,70){\circle*{3}}
\put(120,90){\circle*{3}}
\put(120,110){\circle*{3}}
\put(120,130){\circle*{3}}
\put(35,47){$\times$}
\put(55,47){$\times$}
\put(75,47){$\times$}
\put(95,47){$\times$}
\put(115,47){$\times$}
\put(55,67){$\times$}
\put(75,67){$\times$}
\put(95,67){$\times$}
\put(115,67){$\times$}
\put(75,87){$\times$}
\put(95,87){$\times$}
\put(115,87){$\times$}
\put(95,107){$\times$}
\put(115,107){$\times$}
\put(115,127){$\times$}
\end{picture}
\end{center}
\caption{\label{fi:penrosekeyplus} 
\protect{\footnotesize The image of $\IZ_+ \oplus \IZ_+$ under $T$.   }}
\end{figure}

\begin{figure}[htb]
\begin{center}
\begin{picture}(600,220)(-225,-40)
\put(-180,50){\vector(1,0){130}}
\put(-170,40){\vector(0,1){130}}
\put(-170,50){\vector(1,-1){100}}
\put(30,50){\vector(1,0){130}}
\put(40,40){\vector(0,1){130}}
\put(40,50){\vector(1,1){100}}
\put(-70,40){$a_{n}$}
\put(-183,150){$b_n$}
\put(145,40){$a_{n+1}$}
\put(18,150){$b_{n+1}$}
\put(-175,47){$\times$}
\put(-175,67){$\times$}
\put(-175,87){$\times$}
\put(-175,107){$\times$}
\put(-175,127){$\times$}
\put(-155,47){$\times$}
\put(-155,67){$\times$}
\put(-155,87){$\times$}
\put(-155,107){$\times$}
\put(-155,127){$\times$}
\put(-135,47){$\times$}
\put(-135,67){$\times$}
\put(-135,87){$\times$}
\put(-135,107){$\times$}
\put(-135,127){$\times$}
\put(-115,47){$\times$}
\put(-115,67){$\times$}
\put(-115,87){$\times$}
\put(-115,107){$\times$}
\put(-115,127){$\times$}
\put(-95,47){$\times$}
\put(-95,67){$\times$}
\put(-95,87){$\times$}
\put(-95,107){$\times$}
\put(-95,127){$\times$}
\put(-170,50){\circle*{3}}
\put(-150,50){\circle*{3}}
\put(-130,50){\circle*{3}}
\put(-110,50){\circle*{3}}
\put(-90,50){\circle*{3}}
\put(-150,30){\circle*{3}}
\put(-130,30){\circle*{3}}
\put(-110,30){\circle*{3}}
\put(-90,30){\circle*{3}}
\put(-130,10){\circle*{3}}
\put(-110,10){\circle*{3}}
\put(-90,10){\circle*{3}}
\put(-90,-10){\circle*{3}}
\put(-110,-10){\circle*{3}}
\put(-90,-30){\circle*{3}}
\put(-20,100){$\stackrel{T^{-1}}{ \longleftarrow}$}
\put(40,50){\circle*{3}}
\put(40,70){\circle*{3}}
\put(40,90){\circle*{3}}
\put(40,110){\circle*{3}}
\put(40,130){\circle*{3}}
\put(60,50){\circle*{3}}
\put(60,70){\circle*{3}}
\put(60,90){\circle*{3}}
\put(60,110){\circle*{3}}
\put(60,130){\circle*{3}}
\put(80,50){\circle*{3}}
\put(80,70){\circle*{3}}
\put(80,90){\circle*{3}}
\put(80,110){\circle*{3}}
\put(80,130){\circle*{3}}
\put(100,50){\circle*{3}}
\put(100,70){\circle*{3}}
\put(100,90){\circle*{3}}
\put(100,110){\circle*{3}}
\put(100,130){\circle*{3}}
\put(120,50){\circle*{3}}
\put(120,70){\circle*{3}}
\put(120,90){\circle*{3}}
\put(120,110){\circle*{3}}
\put(120,130){\circle*{3}}
\put(35,47){$\times$}
\put(55,47){$\times$}
\put(75,47){$\times$}
\put(95,47){$\times$}
\put(115,47){$\times$}
\put(55,67){$\times$}
\put(75,67){$\times$}
\put(95,67){$\times$}
\put(115,67){$\times$}
\put(75,87){$\times$}
\put(95,87){$\times$}
\put(115,87){$\times$}
\put(95,107){$\times$}
\put(115,107){$\times$}
\put(115,127){$\times$}
\end{picture}
\end{center}
\caption{\label{fi:penrosekeyminus} 
\protect{\footnotesize The image of $\IZ_+ \oplus \IZ_+$ under $T^{-1}$.   }}
\end{figure}
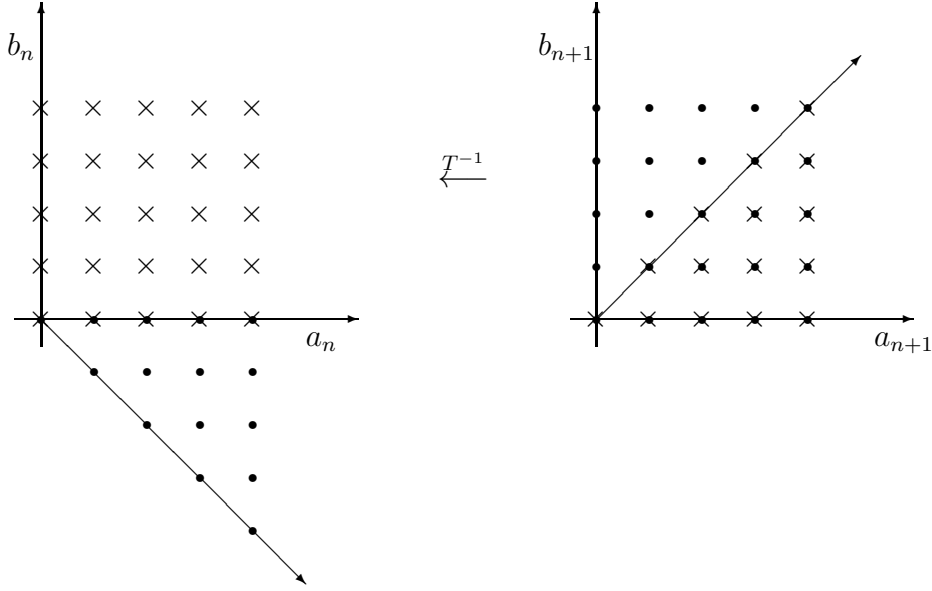\bigskip

\begin{figure}[htb]
\begin{center}
\begin{picture}(600,220)(-350,-40)
\put(-180,50){\vector(1,0){130}}
\put(-170,40){\vector(0,1){130}}
\put(-170,50){\vector(1,-1){100}}
\put(-170,50){\vector(-1,2){55}}
\put(-175,47){$\times$}
\put(-175,67){$\times$}
\put(-175,87){$\times$}
\put(-175,107){$\times$}
\put(-175,127){$\times$}
\put(-155,47){$\times$}
\put(-155,67){$\times$}
\put(-155,87){$\times$}
\put(-155,107){$\times$}
\put(-155,127){$\times$}
\put(-135,47){$\times$}
\put(-135,67){$\times$}
\put(-135,87){$\times$}
\put(-135,107){$\times$}
\put(-135,127){$\times$}
\put(-115,47){$\times$}
\put(-115,67){$\times$}
\put(-115,87){$\times$}
\put(-115,107){$\times$}
\put(-115,127){$\times$}
\put(-95,47){$\times$}
\put(-95,67){$\times$}
\put(-95,87){$\times$}
\put(-95,107){$\times$}
\put(-95,127){$\times$}
\put(-170,50){\circle*{3}}
\put(-170,70){\circle*{3}}
\put(-170,90){\circle*{3}}
\put(-170,110){\circle*{3}}
\put(-170,130){\circle*{3}}
\put(-155,27){$\times$}
\put(-135,27){$\times$}
\put(-115,27){$\times$}
\put(-95,27){$\times$}
\put(-135,7){$\times$}
\put(-115,7){$\times$}
\put(-95,7){$\times$}
\put(-95,-13){$\times$}
\put(-115,-13){$\times$}
\put(-95,-33){$\times$}
\put(-190,90){\circle*{3}}
\put(-190,110){\circle*{3}}
\put(-190,130){\circle*{3}}
\put(-210,130){\circle*{3}}
\end{picture}
\end{center}
\caption{\label{fi:penrosekeyminus2} 
\protect{\footnotesize The image of $\IZ_+ \oplus \IZ_+$ under $T^{-2}$.   }}
\end{figure}
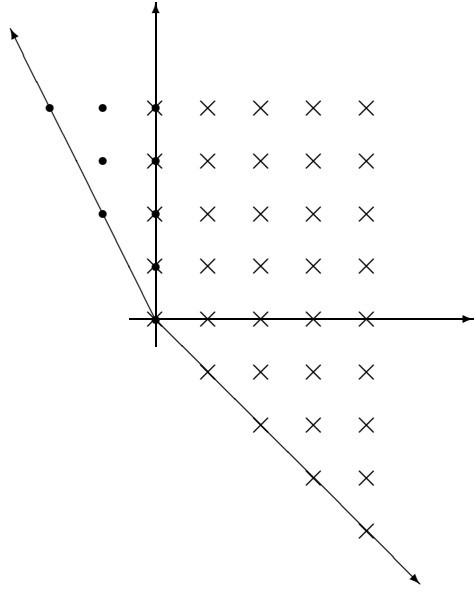\bigskip

\begin{figure}[htb]
\begin{center}
\begin{picture}(600,270)(-350,-90)
\put(-180,50){\vector(1,0){130}}
\put(-170,40){\vector(0,1){130}}
\put(-170,50){\line(3,-5){90}}
\put(-170,50){\line(-3,5){60}}
\put(-175,47){$\times$}
\put(-175,67){$\times$}
\put(-175,87){$\times$}
\put(-175,107){$\times$}
\put(-175,127){$\times$}
\put(-155,47){$\times$}
\put(-155,67){$\times$}
\put(-155,87){$\times$}
\put(-155,107){$\times$}
\put(-155,127){$\times$}
\put(-135,47){$\times$}
\put(-135,67){$\times$}
\put(-135,87){$\times$}
\put(-135,107){$\times$}
\put(-135,127){$\times$}
\put(-115,47){$\times$}
\put(-115,67){$\times$}
\put(-115,87){$\times$}
\put(-115,107){$\times$}
\put(-115,127){$\times$}
\put(-95,47){$\times$}
\put(-95,67){$\times$}
\put(-95,87){$\times$}
\put(-95,107){$\times$}
\put(-95,127){$\times$}
\put(-155,27){$\times$}
\put(-135,27){$\times$}
\put(-115,27){$\times$}
\put(-135,7){$\times$}
\put(-115,7){$\times$}
\put(-115,-13){$\times$}
\put(-135,-13){$\times$}
\put(-115,-33){$\times$}
\put(-115,-53){$\times$}
\put(-95,27){$\times$}
\put(-95,7){$\times$}
\put(-95,-13){$\times$}
\put(-95,-33){$\times$}
\put(-95,-53){$\times$}
\put(-95,-73){$\times$}
\put(-75,-100){${1+\sqrt{5} \over 2} a + b = 0$}
\put(-195,87){$\times$}
\put(-195,107){$\times$}
\put(-195,127){$\times$}
\put(-215,127){$\times$}
\end{picture}
\end{center}
\caption{\label{fi:penkeyzerplu} 
\protect{\footnotesize 
$K_{0+}(\ca)$ for the Penrose tiling. }} 
\end{figure}\bigskip

We shall now evaluate the $K$-theory of the posets associated with the
corresponding Bratteli diagrams in Section \ref{se:bdp}. The strategy will be
the same as the one used for the algebra of the Penrose Tiling and will
consist essentially of three steps:
\begin{itemize}
\item[1.]
Construct, as in Fig.~\ref{fi:ptmatpic} the inclusion maps $T$ in
(\ref{keyinc}) 
from the stable part of the Bratteli diagram.
\item[2.]
Prove that $T$ is invertible over the integer numbers. As a consequence, $K_0$
will be the direct sum of as many copies of $\IZ$ as the number of points
$k_{n_0}$ in the stable level $n_0$ of the corresponding Bratteli diagram of
the poset.
\item[3.]
Identify the subset
$T(K_{0+}(\ca_{n_0})) \subset K_{0+}(\ca_{n_0+1})$ with $K_{0+}(\ca_{n_0})$;

evaluate $T^{-1}(K_{0+}(\ca_{n_0}))$;

get the limit $K_{0+}(\ca) = \lim_{m \ra \infty}
T^{-m}K_{0+}(\ca_{n_0})$.
\end{itemize}

\bexam
The $K$-theory of the poset $\bigvee$. {}From the stable part of the
corresponding Bratteli diagram in Fig.~\ref{fi:veealg}, we get for the
inclusions (\ref{keyinc}), (\ref{keyinc+}) and their inverses the
integer valued 
matrices
\be
T =
\left[
\begin{array}{ccc}
1 & 0 & 0 \\
1 & 1 & 1 \\
0 & 0 & 1
\end{array}
\right]~,~~~
T^{-1} =
\left[
\begin{array}{rrr}
1 & 0 & 0 \\
-1 & 1 & -1 \\
0 & 0 & 1
\end{array}
\right]~.  \label{bigmat}
\ee
Since $T$ is invertible over $\IZ$, from definition (\ref{key}), it follows
that
\be
K_0(\ca(\vee)) = \IZ^3~.
\ee
On the other side, with $(a, b, c) \in \IZ_+^3$, one finds that
\be
T^{-m}
\left(
\begin{array}{l}
a \\
b \\
c
\end{array}
\right) =
\left(
\begin{array}{l}
a \\
b - m(a + c) \\
c
\end{array}
\right)~.
\ee
While $a, c \geq 0$, $b -ma - mc$ can become any negative integer provided
that $a, c \not= 0$. Therefore
\be
\begin{array}{llll}
K_{0+}(\ca(\bigvee)) = \{ (a, b, c) \in \IZ^3
~~|~  & a \in \IZ_+~, c \in \IZ_+ & & \\
     ~                    & b \in \IZ & if &  (a,c) \not= (0,0) \\
     ~                    & b \in \IZ_+ & if & (a,c) = (0,0) ~\}~.
\end{array}
\label{veekey+}
\ee
\eexam
\bexam
The $K$-theory of the poset $\pinco$. {}From the stable part of the
corresponding Bratteli diagram in Fig.~\ref{fi:pinbra}, we get for the
inclusions (\ref{keyinc}), (\ref{keyinc+}) and their inverses the
integer valued 
matrices
\be
T =
\left[
\begin{array}{cccc}
1 & 0 & 0 & 0 \\
1 & 1 & 0 & 1 \\
0 & 0 & 1 & 1 \\
0 & 0 & 0 & 1
\end{array}
\right]~,~~~
T^{-1} =
\left[
\begin{array}{rrrr}
1 & 0 & 0 & 0 \\
-1 & 1 & 0 & -1 \\
0 & 0 & 1 & -1 \\
0 & 0 & 0 & 1
\end{array}
\right]~.   \label{pinmat}
\ee
Since $T$ is invertible over $\IZ$, from definition (\ref{key}), it follows
that
\be
K_0(\ca(\pinco)) = \IZ^4~.
\ee
On the other side, with $(a, b, c, d) \in \IZ_+^4$, one finds that
\be
T^{-m}
\left(
\begin{array}{l}
a \\
b \\
c \\
d
\end{array}
\right) =
\left(
\begin{array}{l}
a \\
b - m(a + d) \\
c - md \\
d
\end{array}
\right)~.
\ee
As a consequence,
\be
\begin{array}{llll}
K_{0+}(\ca(\pinco)) = \{(a, b, c, d) \in \IZ^4
~~|~ & a \in \IZ_+~, d \in \IZ_+ & &  \\
     ~       & b \in \IZ~, c \in \IZ       & if & a \not= 0~, d \not= 0 \\
     ~       & b \in \IZ~, c \in \IZ_+     & if & a \not= 0~, d = 0 \\
     ~       & b \in \IZ_+~, c \in \IZ_+   & if & (a,c) = (0,0)~ \}~.
\end{array}  \label{pinkey+}
\ee
\eexam
\bexam
The $K$-theory of the poset $P_4(S^1)$. {}From the stable part of the
corresponding Bratteli diagram in Fig.~\ref{fi:cirbra}, we get for the
inclusions (\ref{keyinc}), (\ref{keyinc+}) and their inverses the
integer valued 
matrices
\be
T =
\left[
\begin{array}{cccc}
1 & 0 & 0 & 0 \\
1 & 1 & 0 & 1 \\
1 & 0 & 1 & 1 \\
0 & 0 & 0 & 1
\end{array}
\right]~,~~~
T^{-1} =
\left[
\begin{array}{rrrr}
1 & 0 & 0 & 0 \\
-1 & 1 & 0 & -1 \\
-1 & 0 & 1 & -1 \\
0 & 0 & 0 & 1
\end{array}
\right]~.  \label{cirmat}
\ee
Since $T$ is invertible over $\IZ$, from definition (\ref{key}), it follows
that
\be
K_0(\ca(P_4(S^1))) = \IZ^4~.
\ee
On the other side, with $(a, b, c, d) \in \IZ_+^4$, one finds that
\be
T^{-m}
\left(
\begin{array}{l}
a \\
b \\
c \\
d
\end{array}
\right) =
\left(
\begin{array}{l}
a \\
b - m(a + d) \\
c - m(a + d) \\
d
\end{array}
\right)~.
\ee
As a consequence,
\be
\begin{array}{llll}
K_{0+}(\ca(P_4(S^1))) = \{ (a, b, c, d) \in \IZ^4
~~|~ & a \in \IZ_+~, d \in \IZ_+ & &  \\
     ~       & b \in \IZ~, c \in \IZ       & if & a \not= 0 ~~or~~ d \not= 0 \\
     ~       & b \in \IZ_+~, c \in \IZ_+   & if & (a,c) = (0,0)~ \}~.
\end{array} \label{cirkey+}
\ee
Notice that as abelian group, $K_0(\ca(\pinco)) = K_0(\ca(P_4(S^1))) =
\IZ^4$. But as abelian ordered group $(K_0, K_{0+})\ca(\pinco) \not=
(K_0, K_{0+})(\ca(P_4(S^1)))$, since the positive cone $K_{0+}$ is not the
same in the two groups.
\eexam
\bexam
The $K$-theory of the poset $P_6(S^2)$. {}From the stable part of the
corresponding Bratteli diagram in Fig.~\ref{fi:sphbra}, we get for the
inclusions (\ref{keyinc}), (\ref{keyinc+}) and their inverses the
integer valued 
matrices
\be
T =
\left[
\begin{array}{cccccc}
1 & 0 & 0 & 0 & 0 & 0 \\
1 & 1 & 0 & 0 & 0 & 1 \\
1 & 1 & 1 & 0 & 1 & 1 \\
1 & 1 & 0 & 1 & 1 & 1 \\
1 & 0 & 0 & 0 & 1 & 1 \\
0 & 0 & 0 & 0 & 0 & 1
\end{array}
\right]~,~~~
T^{-1} =
\left[
\begin{array}{rrrrrr}
1 & 0 & 0 & 0 & 0 & 0 \\
-1 & 1 & 0 & 0 & 0 & -1 \\
1 & -1 & 1 & 0 & -1 & 1 \\
1 & -1 & 0 & 1 & -1 & 1 \\
-1 & 0 & 0 & 0 & 1 & -1 \\
0 & 0 & 0 & 0 & 0 & 1
\end{array}
\right]~. \label{sphmat}
\ee
Since $T$ is invertible over $\IZ$, from definition (\ref{key}), it follows
that
\be
K_0(\ca(P_4(S^1))) = \IZ^6~.
\ee
On the other side, with $(a, b, c, d, e, f) \in \IZ_+^6$, one finds that
\be
T^{-m}
\left(
\begin{array}{l}
a \\
b \\
c \\
d \\
e \\
f
\end{array}
\right) =
\left(
\begin{array}{l}
a \\
b - m(a + f) \\
c + m^2(a + f) - m(b + e) \\
d + m^2(a + f) - m(b + e) \\
e - m(a + f) \\
f
\end{array}
\right)~.
\ee
As a consequence,
\be
\begin{array}{lll}
K_{0+}(\ca(P_6(S^2))) = \{ (a, b, c, d, e, f) \in \IZ^6 ~~|~~ 
                                      a \in \IZ_+~, f \in \IZ_+ & & \\
~~~~~~~~~~~~~~~~~~~~~~~~~~~~~~~ c \in \IZ~, d
\in \IZ ~~ b \in
\IZ~,
                   e \in \IZ ~~ & if ~~ a \not=0 ~~or~~ f \not= 0 & \\
~~~~~~~~~~~~~~~~~~~~~~~~~~~~~~~~~~~~~~~~~~~~~~~~~~
b \in \IZ_+~, e \in \IZ_+ ~~ & if ~~ (a,f) = (0,0) & \\
~~~~~~~~~~~~~~~~~~~~~~~~~~~~~~~ c \in \IZ_+~, d \in \IZ_+ ~~ & if ~~
(a, b, e, f) = (0, 0, 0, 0)~ \}~. &
\end{array} \\
~ \label{sphkey+}
\ee
\eexam

\sxn{Final Remarks}

As mentioned before, the construction of the $K$-theory groups for
noncommutative lattices is a preliminary step in the classification and
construction of bundles over them and for the theory of characteristic classes.
Notably, one would like to construct non trivial bundles, like, for
instance, the 
analogue of the monopole bundle  over the lattices approximating the
$2$-dimensional sphere and non trivial `topological charges'.  Work in this
direction is in progress and will be  reported elsewhere
\cite{ELTbundles}.

\bigskip

\bigskip

\noindent
{\large \bf Acknowledgments }

Much of this work has been done while G.L. and P.T-S were at ESI in Vienna.
They would like to thank G. Marmo and P. Michor for the invitation and
all people at the Institute for the warm and friendly atmosphere. We are
grateful to A.P. Balachandran, G. Bimonte, F. Lizzi e G. Sparano for many
fruitful discussions and useful advises. We thank U. Bruzzo and R.
Catenacci and J. Varilly for carefully reading the compuscript.
J. Varilly deserves a special mentions for suggestions which improved an early
version.

We thank the `Istituto Italiano per gli Studi Filosofici' in Napoli for
partial support.
The work of P.T-S. was also supported by the Department of Energy,
U.S.A. under
contract number DE-FG-02-84ER40173. The work of G.L. was
partially supported by the Italian `Ministero dell' Universit\`a e
della Ricerca Scientifica'.

\vfill\eject

\end{document}